\documentclass[12pt]{amsart}
\usepackage{amscd}
\usepackage{amssymb}

\numberwithin{equation}{section}

\newcommand{\lab}[1]{\label{#1}}

\newtheorem{Thm}{Theorem}[section]
\newtheorem{Prop}[Thm]{Proposition}
\newtheorem{Lem}[Thm]{Lemma}
\newtheorem{Cor}[Thm]{Corollary}

\theoremstyle{remark}
\newtheorem{Rem}[Thm]{Remark}
\newtheorem{Exa}[Thm]{Example}

\theoremstyle{definition}
\newtheorem*{Dig}{Digression}
\newtheorem{Def}[Thm]{Definition}

\newcommand{\bbR}{{\mathbb{R}}}
\newcommand{\bbQ}{{\mathbb{Q}}}
\newcommand{\cs}{{\rm CS}}

\newcommand{\de}{\partial}
\newcommand{\bpeta}{{\bar p^*\eta}}
\newcommand{\heta}{\hat\eta}
\newcommand{\pfrac}[1]{\genfrac(){}{}{#1}}
\newcommand{\bx}{\mathbf{x}}
\newcommand{\ba}{\mathbf{a}}
\newcommand{\Tr}{\operatorname{Tr}}
\newcommand{\piott}{\pi^{(123)}}
\newcommand{\wC}{\widetilde{C}}
\newcommand{\wcC}{\widetilde{\mathcal{C}}}
\newcommand{\hwC}{\Hat{\Tilde{C}}}
\newcommand{\hwcC}{\Hat{\Tilde{\mathcal{C}}}}

\newcommand{\ord}{\operatorname{ord}}
\newcommand{\sln}{\operatorname{sln}}
\newcommand{\sgn}{\operatorname{sgn}}
\newcommand{\Ad}{\operatorname{Ad}}

\newcommand\qq{\rm}
\newcommand\cmp[1]{{\qq Commun.\ Math.\ Phys.\ \bf #1}}
\newcommand\jmp[1]{{\qq J.\ Math.\ Phys.\ \bf #1}}
\newcommand\pl[1]{{\qq Phys.\ Lett.\ \bf #1}}

\newcommand\anm[1]{{\qq Ann.\ Math.\ \bf #1}}

\newcommand\jdg[1]{{\qq J.\ Diff.\ Geom.\ \bf #1}}

\begin{document}
\title{Integral Invariants of 3-Manifolds}
\author[R.~Bott]{Raoul Bott}
\address{Department of Mathematics\\
Harvard University\\
Cambridge, MA~02138, USA}
\email{bott@math.harvard.edu}

\author[A.~S.~Cattaneo]{Alberto S.~Cattaneo}
\address{
Lyman Laboratory of Physics\\
Harvard University\\
Cambridge, MA~02138, USA}
\email{cattaneo@math.harvard.edu}

\date{July 2, 1997}

\thanks{We thankfully acknowledge support from NSF for R.B. and
from INFN (grant No.\ 5565/95) 
and from DOE (grant No.\ DE-FG02-94ER25228,
Amendment No.\ A003) for A.S.C}

\begin{abstract}
This note describes an invariant of rational homology 3-spheres
in terms of configuration space integrals which in some sense
lies between the invariants of Axelrod and Singer \cite{AS}
and those of Kontsevich \cite{K}.
\end{abstract}
\maketitle

\section{Introduction}
In their seminal paper of 1994, \cite{AS}, Axelrod and Singer 
showed that the asymptotics of the Chern--Simons theory led to
a series of $C^\infty$-invariants associated to triples
$\{M;f;\rho\}$ with $M$ a smooth homology 3-sphere,
$f$ a homotopy class of framings of $M$, and $\rho$ an
``acyclic'' conjugacy class of orthogonal representations
of $\pi_1(M)$. That is, the cohomology $H^*(M;\Ad\rho)$ of $M$
relative to the local system associated to $\Ad\rho$ vanishes.

The primary purpose of this note is to show that the basic ideas
of their paper can be adapted quite easily---but not quite
trivially---to yield invariants of smooth, framed 3-dimensional
homology spheres as such. Put differently, we will present a
treatment somewhat analogous to theirs for the trivial
representation of $\pi_1(M)$.
We say somewhat because in our work we have put aside the
physics inspired aspects of Axelrod and Singer's paper.
Instead we have simply taken our task to be the production
of invariants of framed manifolds $(M;f)$ out of some fixed
Riemannian structure on $M$.

There is of course Kontsevich's solution by ``softer methods"
to the problem of finding the residual invariants of the
Chern--Simons theory at the trivial representation.
In a note (see \cite{K}), Kontsevich sketched how to define
a series of invariants for ``framed'' 3-dimensional
homology spheres, and developed his ``graph cohomology''
to explain the combinatorial diversity of these invariants.

In 1995 Taubes \cite{CT} carefully investigated the first of
these Kontsevich invariants---the one associated to
the $\Theta$-graph---and clarified the appropriate concept
of ``framing'' for all of the graph cohomology.
He coined the term ``singular framings'' for them,
and they differ from ordinary framings in that they exhibit
a singularity at one point of $M$.

Our invariants---which are less ``soft" than Kontsevich's in the
sense that they do depend on Riemannian concepts---are therefore,
on the face of it,
different from his. But they are also indexed by cocycles, $\Gamma$,
in an appropriate graph cohomology, and structurally take the form:
\begin{equation}
I_\Gamma(M,f) = A_\Gamma(M) + \phi(\Gamma)\,\cs(M,f).
\lab{defI}
\end{equation}
Here $A_\Gamma(M)$ is an integral over the configuration spaces
specified by $\Gamma$ and defined by a fixed Riemannian structure
on $M$, $\phi(\Gamma)$ is a real number universally associated
to $\Gamma$, and $\cs(M,f)$ denotes the Chern--Simons integral
of the Levi-Civita connection of $M$ relative to the frame $f$.

The Axelrod--Singer invariants for a flat connection exhibit a 
similar dependance on the framings, and \eqref{defI} is also in
general conformity with the self-linking invariants in knot
theory---as described in \cite{BT} (see also \cite{AF}). 
There the invariants
of a knot $K\subset\bbR^3$ are described as
\begin{equation}
I_\Gamma(K) = A_\Gamma(K) + \mu(\Gamma)\cdot\text{self-linking of $K$},
\end{equation}
where $I_\Gamma(K)$ is a configuration space integral which
is corrected by an anomalous term which is a multiple of
the self-linking of $K$.

Similarly we now obtain invariants of oriented homology 3-spheres,
one for every connected cocycle $\Gamma$, of the form
\begin{equation}
J_\Gamma(M) = A_\Gamma(M) - 4\,\phi(\Gamma)\, A_\Theta(M),
\end{equation}
so that $A_\Theta(M)$ is seen to play the role of the self-linking
integral in knot theory.

Although the invariants of \cite{AS} and \cite{BT}, as well as
the ones described here, are all spin-offs from Witten's \cite{W}
original Chern--Simons invariants for homology 3-spheres,
it seems to us that, from a purely mathematical point of view,
they have now, in retrospect, even older antecedents.
These are the ``iterated integrals" of Chen, or---even
older---the Adams constructions for the loop-space of a space.

Quite generally, the principle of these constructions is to
describe the cohomology of a function-space $F={\rm Map}(X,Y)$
in terms of the various evaluation maps:
\[
{\rm Map}(X,Y) \times X^n \to Y^n.
\]
When we are dealing with corresponding spaces of imbeddings,
or diffeomorphisms, then the configuration spaces enter the
discussion quite naturally, and give rise to new invariants of
the type we have been discussing.

In this context it is also possible to extend our considerations
to local systems on $M$, and derive similar invariants, all
governed by some graph cohomology.
{}From this point of view the original invariants of
Axelrod--Singer are associated to the ``Feynman cocycles"
of this cohomology.

We hope to explore these ideas in a subsequent 
paper \cite{BC2}. Here
we will only deal with the constant coefficient case
and the plan of this note is as follows:
In section \ref{sec-cc} we review some facts from the theory
of characteristic classes.
In section \ref{sec-inv} we describe the $\Theta$-invariant
explicitly, but implicitly rely on the description of
configuration spaces as developed in \cite{AS}, in analogy
with the corresponding algebraic construction given by
Fulton and MacPherson in \cite{FM}.
In section \ref{sec-hi} we discuss the higher invariants
while the last, fifth, section is devoted to extending
the results of \cite{BT} to knots in general homology
3-spheres.

\subsection*{Acknowledgements}
We are indebted for very useful conversations pertaining to
these matters with Scott Axelrod, Robin Forman, Stavros
Garoufalidis and Cliff Taubes.

\section{Review of characteristic classes of $SO(n)$}
\lab{sec-cc}
Consider an oriented vector bundle $E$ with odd fiber dimension,
$n=2k+1$, over a base space $M$. Also let $S(E)$ denote the associated
sphere bundle to $E$, which we may consider to be the space of rays
in $E$; or, if $E$ is given a Riemannian structure, as the unit sphere
bundle of $E$. In any case $S(E)$ has even fiber dimension $2k$ over
$M$, and this together with the orientability of $E$ allows one to
specify a canonical integral generator of the rational cohomology
of $S(E)$ as a module over $H^*(M)$. Namely, we consider the ``tangent
bundle along the fiber,'' $T_FS(E)$, of $S(E)$. This, being an even dimensional
oriented bundle, has a canonical Euler class:
\begin{equation}
e = e(T_FS(E))\in H^{2k}(S(E)),
\end{equation}
which restricts to twice the generator of $H^{2k}(S^{2k})$ on each fiber,
because the Euler number of $S^{2k}$ is 2.

But then it follows from general principles that $e$ generates
$H^*(S(E))$ over $H^*(M)$ over the rationals.

Concerning the generator $e$ we have the following lemma, which
in some sense explains the Chern--Simons term in our subsequent
construction.
\begin{Lem}
Let $\pi_*$ denote integration along the fiber in the bundle
$S(E)$ over $M$. Then
\begin{equation}
\pi_* e^3 = 2\,p_k(E),
\lab{pie}
\end{equation}
where $p_k$ denotes the $k$-th Pontrjagin class of $E$.
In fact one has, quite generally:
\begin{align}
\pi_* e^{2s+1} &= 2\, (p_k(E))^s,\\
\intertext{and}
\pi_* e^{2s} &=0,
\end{align}
for $s=1,2,\dotsc$.
\lab{lem-pie}
\end{Lem}
\begin{proof}
It suffices to prove these formulae for the universal sphere bundle
over the classifying space $BSO(n)$, $n=2k+1$, that is for the
fibering
\begin{equation}
\begin{CD}
BSO(2k)\\
@V{\pi}VV\\
BSO(2k+1)
\end{CD}
\end{equation}
with fiber $SO(2k+1)/SO(2k)=S^{2k}$.

Here we can keep track of the rational cohomology of the spaces involved
by choosing a maximal torus $T=(S^1)^k$ for $SO(2k+1)$ in the usual
manner, so that $T$ corresponds to diagonal $2\times2$ blocks ending
with a 1 in the last diagonal position. In this way $H^*(BT)$
becomes identified with the polynomial ring $\bbQ[x_1,\dots,x_k]$
and the Weyl group of $SO(2k+1)$ acts on this cohomology by
1) permutations of the $x_i$, and 2) changes of sign $x_i\to\pm x_i$.
On the other hand the Weyl group of $SO(2k)$ acts as the subgroup
which 1) permutes the $x_i$, and 2) allows only even changes of sign
$x_i\to\epsilon_ix_i$, $\epsilon_i=\pm1$, with $\prod\epsilon_i=1$.

It follows that the invariants of $H^*(BT)$ under $SO(2k+1)$ are
given by the invariant polynomials 
$\sigma_r=\sigma_r(x_1^2,\dots,x_k^2)$, $r\le k$, while those
invariant under $SO(2k)$ are generated by $\sigma_r$ and an
additional element
\begin{equation}
e = x_1\dots x_k\in H^{2k}(BT).
\lab{e}
\end{equation}

{}From the well-known identification of $H^*(BSO(2k+1))$ and
$H^*(BSO(2k))$ with these rings of invariants respectively,
we conclude that:
\begin{equation}
e^2 = \pi^*(x_1^2\dots x_k^2) = \pi^*\sigma_k \qquad
\text{in \eqref{e}}.
\lab{e2}
\end{equation}
But then
\begin{equation}
\pi_*e^3 = \pi_*(e\,\pi^*\sigma_k) = 2\,\sigma_k
\end{equation}
by the permanence relation and the fact that $\pi_*e=2$ 
remarked upon
earlier. If we take $\sigma_r$ to be universal Pontrjagin 
class---as opposed to the convention $p_r=(-1)^r\,\sigma_r$---then
\eqref{e2} implies \eqref{pie}, and the general case follows
similarly from $e^{2s+1}=e\,\pi^*(\sigma_k)^s$.
\end{proof}

\section{The simplest invariant}
\lab{sec-inv}
From now on we will only consider a 3-dimensional 
rational homology sphere $M$.
The boundary of the configuration space of two points in
$M$, $C_2(M)$, is then isomorphic to
the 2-sphere bundle $S(TM)$ over $M$. In the previous
section we have seen how to construct a vertical generator $e$
in a sphere bundle using Riemannian concepts. 
In this section we will give this generator explicitly
as an element of $\Omega^2(\de C_2(M))$. 
In the de Rham theory we can divide by 2 and so we will actually
describe $\eta=e/2$.
Then we will extend
it to the whole of $\Omega^2(C_2(M))$ and will show that its
differential is the Poincar\'e dual of the diagonal in $M\times M$.
The next step will be to use this element to construct a 
closed form in
$\Omega^2(C_3(M))$. It is precisely through this form that we will be
able to write the simplest invariant of the rational homology sphere
$M$ as
an integral over $C_3(M)$ (notice that in \cite{AS}, 
\cite{K} and \cite{CT} the ``corresponding'' $\Theta$-invariant is written
as an integral over $C_2(M)$). Finally, we will prove that, apart
from an anomalous term which we compute explicitly, this is 
actually an invariant.

\subsection{The generator of $H^2(\de C_2(M))$}
We may consider $\de C_2(M)$ as the sphere bundle
$P\times_{SO(3)}S^2\to M$ where $P\to M$ is the orthonormal frame bundle
of $TM$ with respect to some fixed Riemannian structure, so that 
we have the commutative diagram:
\[
\begin{CD}
P\times_{SO(3)}S^2 @<{\bar p}<< P\times S^2\\
@V{\pi^\de}VV @VV{\bar\pi^\de}V\\
M @<<{p}< P
\end{CD}
\]
Note that here $\bar\pi^\de$ is a morphism of principal 
$SO(3)$-bundles while $\bar p$ is a corresponding morphism
of $S^2$-bundles.

We will write our class $\eta$, or more precisely $\bpeta$, as
a closed form in $\Omega^2_{\text{basic}}(P\times S^2)$ such that
$\bar\pi^\de_*\bpeta = 1$. For the $(0,2)$ component of $\bpeta$
we choose the $SO(3)$-invariant volume element
\[
\omega = x\, dy\, dz + y\, dz\, dx + z\, dx\, dy =
\frac12\epsilon^{ijk}\,x_i\,dx_j\,dx_k,
\]
which satisfies $\int_{S^2}\omega=4\pi$.
The $SO(3)$-action on $S^2$ is given by the vector fields 
\[
X_i = {\epsilon_{ij}}^k\, x_j \frac{\de}{\de x_k}.
\]
We have
\[
L_{X_i}\,\omega = 0, \quad \iota_{X_i}\,\omega = d x_i.
\]
Let $\{\xi_1,\xi_2,\xi_3\}$ be the basis of $\mathfrak{so}(3)$
that corresponds to the vector fields $X_1,X_2,X_3$; that is,
\[
(\xi_i)_{jk} = \epsilon_{ijk}.
\]
Then a connection $\theta$ on $P$ can be expanded in this basis
as
\[
\theta = \theta^i\xi_i,
\]
so that, by the definition of a connection, we have
\[
L_{X_i}\,\theta^j = -\epsilon_{ijk}\, \theta^k,\quad
\iota_{X_i}\,\theta^j = \delta_i^j.
\]

In the following we will write $\theta^i$ also for the pullback
of $\theta^i$ to $P\times S^2$, and similarly consider the coordinates
$x_i$ of $\bbR^3$ also as pulled back to $P\times S^2$.
This understood consider
the invariant 1-form $\theta^i x_i$. It then follows from this
invariance that
\[
\iota_{X_i}\, d(\theta^jx_j) = -d(\iota_{X_i}\,\theta^jx_j) = -dx_i.
\]
Next, we have the following
\begin{Prop}
The 2-form
\begin{equation}
\bar\eta =\dfrac{\omega + d(\theta^ix_i)}{4\pi},
\quad\quad x^ix_i=1,
\lab{bpeta}
\end{equation}
is basic in $P\times S^2\to P\times_{SO(3)}S^2$. Moreover,
if we write $\bar\eta=\bpeta$, then
\begin{equation}
\pi^\de_*\eta = 1,\quad \pi^\de_*\eta^2 = 0.
\lab{bpeta10}
\end{equation}
Finally, if $\phi$ is the automorphism of the bundle
$P\times_{SO(3)}S^2\to M$ given by the antipodal map
on the fiber, then
\begin{equation}
\phi^*\eta = -\eta.
\lab{bpetaphi}
\end{equation}
\lab{thm-eta}
\end{Prop}
\begin{proof}
Eqn.\ \eqref{bpeta} follows directly from the previous discussion.

To prove \eqref{bpeta10}, we
notice that $p^*\pi^\de_*=\bar\pi^\de_*\bar p^*$. Thus,
the first identity is just a consequence of the fact
that
\[
\bar\pi^\de_*[\omega+ d(\theta^ix_i)]
=\bar\pi^\de_*\omega=4\pi.
\] 
For the second identity we compute
\begin{multline*}
(4\pi)^2\bar\pi^\de_*(\bpeta)^2 =
\bar\pi^\de_*[2\omega d(\theta^ix_i) +
d(\theta^i x_i)\,d(\theta^jx_j)]=\\
2d\theta^i\bar\pi^\de_*(\omega x_i)-
\theta^i\theta^j\bar\pi^\de_*(dx_idx_j)=0.
\end{multline*}
For the last identity, notice that the integral of $\omega x_i$ 
vanishes by symmetry and that $dx_idx_j$ is exact.

Finally, to prove \eqref{bpetaphi}, we consider the automorphism
$\bar\phi$ of $P\times S^2 \to P$ obtained by the antipodal map
on $S^2$:
\[
\bar\phi x_i = -x_i.
\] 
Clearly, $\bar\phi^*\bpeta= -\bpeta$. Moreover,
$\phi\bar p = \bar p\bar\phi$. 
Then
\[
\bar p^*\phi^*\eta=\bar\phi^*\bpeta=-\bpeta,
\]
which implies \eqref{bpetaphi}.
\end{proof}

\begin{Dig}
If we introduce the covariant derivative $D$,
\[
(D x)_i = dx_i-{\epsilon_{ij}}^k \theta^jx_k,
\]
and the curvature $F = F^i\xi_i$,
\begin{equation}
F^i = d\theta^i-\frac12{\epsilon^i}_{jk}\theta^j\theta^k,
\lab{defF}
\end{equation}
then a straightforward but tedious
computation shows that we can also write
\begin{equation}
\bpeta = \frac{\frac12\epsilon^{ijk}\,x_i(Dx)_j(Dx)_k 
+ F^ix_i}
{4\pi}.
\lab{bpeta'}
\tag{\ref{bpeta}$'$}
\end{equation}
In fact,
\begin{multline*}
\frac12\epsilon^{ijk}\,x_i(Dx)_j(Dx)_k =
\frac12\epsilon^{ijk}\,x_i\,dx_j\,dx_k +\\
-\frac12\epsilon^{ijk}\,x_i\,{\epsilon_{jr}}^s\,\theta^r\,x_s\,dx_k
-\frac12\epsilon^{ijk}\,x_i\,dx_j\,{\epsilon_{kr}}^s\,\theta^r\,x_s+\\
+\frac12\epsilon^{ijk}\,x_i\,{\epsilon_{jr}}^s\,\theta^r\,x_s\,
{\epsilon_{km}}^n\,\theta^m\,x_n=
\omega -\theta^i\,dx_i 
+\frac12{\epsilon^i}_{jk}\,x_i\,\theta^j\,\theta^k.
\end{multline*}
To obtain the last identity we have used
\[
\sum_i\epsilon_{ijk}\,\epsilon_{irs} = \delta_{jr}\delta_{ks}
-\delta_{js}\delta_{kr},
\]
and the constraint $x^ix_i=1$.

Note also that, if
one considers $\bx = x^i\xi_i$ as an element of the fundamental 
representation of $\mathfrak{so}(3)$, then, by using the identities
\[
\Tr \xi_i\xi_j = -2\delta_{ij},\quad
\Tr \xi_i[\xi_j,\xi_k] = 2\epsilon_{ijk},
\]
one can rewrite \eqref{bpeta} and \eqref{bpeta'} as
\begin{equation}
\bpeta=\frac{\omega-\frac12 d\Tr(\theta\,\bx)}{4\pi}=
\frac{\Tr(\bx\, D\bx\, D\bx - F\,\bx)}{8\pi}.
\lab{bpetaTr}
\end{equation}
End of the digression.
\end{Dig}

\subsection{The extension to $C_2(M)$}
First we want extend our form $\eta$ to a small neighborhood
$U$ of $\de C_2(M)$. We can think of this neighborhood as
the complement $TM'$ of the zero section of the tangent bundle $TM$.
We still have an $SO(3)$-bundle
\[
P\times(\bbR\backslash0)\xrightarrow{\bar p} TM'.
\]
Scaling each coordinate $x_i$ in $\bbR\backslash0$
by $r = (x_1)^2+(x_2)^2+(x_3)^2$, we get the closed, basic form
\begin{equation}
\bpeta =\frac\omega{4\pi r^3} + \frac1{4\pi}
d\pfrac{\theta^ix_i}r,
\end{equation}
with $\omega$ defined as before.
Then we consider a second neighborhood $V$ containing $U$ and contained
in $C_2(M)$, and choose a function $\rho$ on $C_2(M)$ that is constant
and equal to $-1$ in $U$ and constant and equal to $0$ in the complement
of $V$. It is then clear that $d(\rho\eta)$ represents a class in 
$H^3(M\times M)=H^3(M)\otimes H^3(M)$. 

Let us denote by $\pi_1$ and $\pi_2$ the
two natural projections from $M\times M$ to $M$, i.e.,
\[
\pi_i(m_1,m_2) = m_i,
\]
and by  $v$ a unit volume form on $M$ ({\em not necessarily the
volume form determined by the metric}). In fact any $v$ with 
$\int_M v =1$ would do, and we will use the term ``unit volume"
form in this sense throughout.
Then the generators
of $H^3(M\times M)$ are $v_1$ and $v_2$, defined by
\[
v_i = \pi_i^* v,
\]
and we can write $[d(\rho\eta)]=c_1v_1+c_2v_2$ for some 
constants $c_1$ and $c_2$.
Since
\[
\int_{C_2(M)} d(\rho\eta) v_i =
\int_{\de C_2(M)} \rho\eta v_i =
-\int_{\de C_2(M)}\eta v_i =
-\int_M v = -1,
\]
we see that actually
\[
[d(\rho\eta)] = v_2-v_1;
\]
that is, $d(\rho\eta)$ represents
the Poincar\'e dual of the diagonal in $M\times M$.
This means that there exists a form $\alpha\in\Omega^2(M\times M)$
such that
\begin{equation}
d(\rho\eta) =v_2-v_1-d\alpha.
\lab{dreta}
\end{equation}
Now consider the involution
\begin{equation}
T : \begin{aligned}[t]
C_2(M) &\to C_2(M)\\
(m_1,m_2) &\mapsto (m_2,m_1),
\end{aligned}
\lab{defT}
\end{equation}
and its analog on $M\times M$, which we still denote by $T$.
Since $T$ restricted to $\de C_2(M)$ is the automorphism $\phi$
considered in Prop.\ \ref{thm-eta}, then
\[
T^*(\rho\eta)=-\rho\eta,
\]
provided we choose $\rho$ symmetric (e.g., we can take $\rho$ to
be a function of the distance between $m_1$ and $m_2$).
It follows that in \eqref{dreta} we can choose $\alpha$
such that
\[
T^*\alpha = -\alpha.
\]
Define
\begin{equation}
\heta = \rho\eta + \alpha \in\Omega^2(C_2(M)).
\lab{defheta}
\end{equation}
We have therefore proved the following
\begin{Prop}
There exist forms $\heta\in\Omega^2(C_2(M))$ with the
following three properties:
\begin{subequations}\lab{propheta}
\begin{gather}
\pi^\de_*\iota_\de^*\heta = -1,\\
d\heta = v_2-v_1,\\
T^*\heta =-\heta.
\end{gather}
\end{subequations}
Moreover, there exist forms $\heta$ with the additional
property
\begin{equation}
\iota_\de^*\heta=-\eta.
\lab{propheta'}
\end{equation}
\end{Prop}
Here $\iota_\de$ is the inclusion $\de C_2(M)\hookrightarrow C_2(M)$.
\begin{Rem}
A metric, a compatible connection and a unit volume
form are not enough to determine a unique $\heta$, for
\begin{equation}
\heta'=\heta+d\beta,
\lab{heta'}
\end{equation}
with $\beta\in\Omega^1(C_2(M))$
such that $T^*\beta=-\beta$, still satisfies \eqref{propheta}.
If we moreover want $\heta$ to satisfy \eqref{propheta'},
then we must also put the restriction that $\iota_\de^*\beta=0$.
\end{Rem}

\begin{Dig}[The Riemannian parametrix]
Given a Riemannian structure $g$ on a manifold $M$, 
a linear operator
\[
P_g : \Omega^*(M) \to \Omega^{*-1}(M)
\]
with the property that
\begin{equation}
dP_g + P_g d = 1-\pi_h,
\lab{defPg}
\end{equation}
where $\pi_h$ is the orthogonal projection onto the harmonic forms,
will be called a Riemannian parametrix.
Of course \eqref{defPg} does not define a unique $P_g$, for
\begin{equation}
P_g'=P_g + dQ-Qd
\lab{Pg'}
\end{equation}
still satisfies it for any $Q:\Omega^*(M) \to \Omega^{*-2}(M)$. 

The harmonic projection can be written as a convolution on
$M\times M$ (or $C_2(M)$) as
\[
\pi_h\alpha = \pi_{2*}(\eta_\Delta\,\pi_1^*\alpha)
\]
where $\eta_\Delta$ is the representative of the Poincar\'e dual of the diagonal in $M\times M$ determined by the metric on $M$. In the case
when $M$ is a rational homology sphere we have $\eta_\Delta=v_2-v_1$.

Now we have the following
\begin{Prop}
A form $\heta\in\Omega^2(C_2(M))$ satisfying \eqref{propheta} with
$v$ the volume form determined by the metric $g$ 
is the Schwartz kernel for a Riemannian parametrix $P_g$. 
More precisely, given a form $\alpha\in\Omega^*(M)$, the operator 
$P_g$ defined by
\begin{equation}
P_g\alpha = -\pi_{2*}(\heta\,\pi_1^*\alpha)
\lab{Pgheta}
\end{equation}
satisfies \eqref{defPg}.\lab{prop-Pgheta}
\end{Prop}
\begin{proof}
We need the following generalization of Stokes' formula:
\begin{equation}
d\pi_{2*} = -\pi_{2*}d+\pi^\de_*\iota_\de^*,
\lab{dpi*}
\end{equation}
which holds in the case of an odd-dimensional fiber with boundary.
It follows that
\begin{multline*}
dP_g\alpha+P_gd\alpha =
\pi_{2*}(d\heta\,\pi_1^*\alpha)-
\pi_*^\de\iota_\de^*(\heta\,\pi_1^*\alpha)=\\
-v\,\pi_{2*}\pi_1^*\alpha-
\pi_{2*}\pi_1^*(v\,\alpha)+
\alpha\,\pi_*^\de\eta=
(1-\pi_h)\alpha,
\end{multline*}
where we have also used the fact that
$\iota^*_\de\pi_i^*=\pi^{\de*}$.
\end{proof}
\begin{Rem}
To define the Riemannian parametrix, we have only used 
properties \eqref{propheta}. The additional property
\eqref{propheta'} will be crucial to define the manifold
invariants. Notice, moreover, that the freedom \eqref{heta'}
in defining $\heta$ corresponds to the freedom \eqref{Pg'}
in defining $P_g$.
\end{Rem}
\begin{Rem}
A particular choice of $P_g$ is given by
$\bar P_g = d^*\circ G$,
where $G$ is the inverse of $\square + \pi_h$ and $\square$ is the
Laplace operator determined by $g$. We will not concentrate our
interest on this particular Riemannian parametrix---as was the 
case in \cite{AS}---but will stick to the general case.
In \cite{AS} it is precisely
the Schwartz kernel for this Riemannian parametrix
$\bar P_g$ that is constructed, and found to be represented
on the boundary precisely by the form $\eta$ we have been
considering.
Close to the boundary there
are corrections which are continuous but not smooth as forms on
$M\times M$ (corresponding to the singular part of $G$).
These forms, however, become smooth when lifted to $C_2(M)$.
Then, with a suitable choice of $\beta\in\Omega^1(C_2(M))$ in
\eqref{heta'}, we can recover the $\heta$ representing $\bar P_g$.
\end{Rem}
End of the digression.
\end{Dig}

\subsection{Extension to $C_3(M)$} 
Consider the three natural
projections $\pi_1$, $\pi_2$ and $\pi_3$ from $C_3(M)$ to
$M$ given by
\[
\pi_i(m_1,m_2,m_3) = m_i,
\]
and call
\[
v_i = \pi_i^* v.
\]
Then consider the three natural
projections $\pi_{12}$, $\pi_{23}$ and $\pi_{13}$ from $C_3(M)$
to $C_2(M)$:
\[
\pi_{ij}(m_1,m_2,m_3) = (m_i,m_j),\quad
{1\le i<j\le 3},
\]
and define
\[
\pi_{ji} = T\pi_{ij},
\]
where $T$ is the involution defined in \eqref{defT}.
We will denote by
\[
\heta_{ij} =\pi_{ij}^*\heta
\]
the pullbacks of the form $\heta$ defined in \eqref{defheta}.
We can recast the properties of $\heta$ as
\begin{equation}
\begin{split}
d\heta_{ij} &= v_j-v_i,\\
\heta_{ji} &= -\heta_{ij}.
\end{split}
\lab{prophetaij}
\end{equation}
Finally, introduce
\begin{equation}
\heta_{ijk} = \heta_{ij}+\heta_{jk}+\heta_{ki},
\lab{hetaijk}
\end{equation}
for $i$, $j$ and $k$ different from each other.
A simple consequence of \eqref{prophetaij} is:
\begin{equation}
\begin{split}
d\heta_{ijk} &=0,\\
\heta_{ijk} &= \epsilon_{ijk}\,\heta_{123}.
\end{split}
\lab{propheta123}
\end{equation}
This way we have constructed a closed form in $\Omega^2(C_3(M))$.

\begin{Rem}
The form $\heta$ depends on the choice of the unit volume
form. In fact, if we pick up
a different volume form $v'=v+d\tau$, then, by (\ref{propheta}b),
we must replace
$\heta_{ij}$ by $\heta_{ij}'=\heta_{ij}+\tau_j-\tau_i$.
By \eqref{hetaijk}, we see that $\heta_{ijk}$ is unchanged. 
\end{Rem}

We have not used property \eqref{propheta'} yet. First notice that
the boundary of $C_3(M)$ has four faces of codimension one, which
we denote by $(12)$, $(23)$, $(31)$ and $(123)$, by indicating the
underlying diagonal. Then it follows that
\begin{subequations}\lab{propheta123de}
\begin{align}
\iota_{(12)}^*\heta_{123} &= -\eta_{12},\\
\iota_{(23)}^*\heta_{123} &= -\eta_{23},\\
\iota_{(31)}^*\heta_{123} &= -\eta_{31},\\
\intertext{and}
\iota_{(123)}^*\heta_{123} &= -(\eta_{12} + \eta_{23} + \eta_{31}).
\end{align}
\end{subequations}
Here by $\eta_{ij}$ we mean the pullback of the form 
$\eta\in\de C_2(M)$ by the restriction to the boundary of the maps
$\pi_{ij}$.

More precisely, a face like (12) is a sphere bundle over $C_2(M)$.
If we denote by $m_1$ the point in $C_2(M)$ where the collapse
has happened, then (12) can be expressed as
$\pi_1^{-1}\de C_2(M)$, where $\pi_1$ is the corresponding projection
$C_2(M)\to M$. Then $\eta_{12}=\pi_1^*\eta$. Similarly for the faces
$(23)$ and $(31)$.

The face $(123)$ is a bundle over $M$ whose fiber $F$ is given by 
$C_3(\bbR^3)$ modulo global translations and scalings. If we
denote by $\bx_1$, $\bx_2$ and $\bx_3$ the coordinates of $F$,
then we have the projections
\begin{equation}
\pi_{ij} : \begin{array}[t]{ccc}
F &\to& S^2,\\
(\bx_1,\bx_2,\bx_3) &\mapsto& \dfrac{\bx_j-\bx_i}{|\bx_j-\bx_i|},
\end{array} \quad\quad i\not=j,
\lab{piF}
\end{equation}
and their trivial extension to $P\times F\to P\times S^2$. Since
they are equivariant, they descend to $P\times_{SO(3)}F\to
P\times_{SO(3)}\times S^2=\de C_2(M)$. Then $\eta_{ij}=\pi_{ij}^*\eta$.

\begin{Rem}
The form $\heta$ is defined up to the differential of
a 1-form that vanishes on the
boundary. Under the transformation \eqref{heta'}, we have
$\heta_{123}' = \heta_{123} + d(\beta_{12}+\beta_{23}+\beta_{31})$.
More generally, since the properties we are interested in are
\eqref{propheta123} and \eqref{propheta123de},
we can allow the addition of any exact term,
\[
\heta_{123}' = \heta_{123} + d\beta,
\]
with $\beta\in\Omega^1(C_3(M))$ and vanishing on the boundary.
\end{Rem}

\subsection{The simplest invariant}
We now have all the necessary elements to define
the configuration space integral:
\begin{equation}
A_\Theta \doteq \int_{C_3(M)} \heta_{123}^3\,v_1.
\lab{ATheta}
\end{equation}
The apparent asymmetry in the choice of $v_1$ can be removed if we
notice that, by cyclically exchanging the three points in $C_3(M)$,
we also have
\[
A_\Theta(M) = \int_{C_3(M)} \heta_{123}^3\,v_2 =
\int_{C_3(M)} \heta_{123}^3\,v_3 =
\frac13\int_{C_3(M)} \heta_{123}^3\,(v_1+v_2+v_3).
\]
The definition of $A_\Theta$ relies on many choices:
a metric, a connection compatible with that metric and a unit volume
form; moreover, $\heta_{123}$ is defined up to the differential
of a 1-form that
vanishes on the boundary. The last freedom is immediately seen
not to have consequences on $A_\Theta$ since $\heta_{123}$ and
$v_3$ are closed forms.
As we will see in the next subsection, $A_\Theta$ is not
completely independent of all the other choices. However, we will
be able to prove the following
\begin{Thm}
Given a section $f$ of the frame bundle $P$, the combination
\begin{equation}
I_\Theta(M,f) = A_\Theta(M) +\frac14\, \cs(M,f),
\lab{ITheta}
\end{equation}
is independent of all the choices involved (except for the framing).
Here
\begin{multline}
\cs(M,f)=-\frac1{8\pi^2}\int_Mf^*\Tr\left({
\theta\,d\theta +\frac23\,\theta^3}\right)=\\
\frac1{4\pi^2}\int_M f^*\left({
\theta^id\theta_i-\frac13\,\epsilon_{ijk}\theta^i\theta^j\theta^k
}\right),
\end{multline}
is the Chern--Simons integral of the same metric
connection used to define $\eta$.

Thus, $I_\Theta(M,f)$ is an invariant for the framed rational homology
sphere $(M,f)$.
\lab{thm-Theta}
\end{Thm}

\begin{Rem}
In an $SO(3)$-bundle, it is {\em half} the Chern--Simons form
that restricted to the fiber yields the integral 
generator \cite{CS}.
Therefore, the Chern--Simons term is defined up to an even
integer, and the $\Theta$-invariant $I_\Theta$ up
to half an integer.
\lab{rem-CS}
\end{Rem}

\begin{Rem}
So far we have considered $v$ to be a unit volume form (not
necessarily determined by the metric). We can drop this
assumption defining the invariant as
\[
I_\Theta(M,f) = \frac1{V^4}\,A_\Theta(M) +\frac14\, \cs(M,f),
\]
where $V=\int_M v$. Notice that, for $\heta$ to satisfy
\eqref{propheta}, we must now take the function $\rho$ in
\eqref{defheta} to be constant and equal to $-V$ close to the
boundary.
\end{Rem}

If one expands $\heta_{123}^3$ in terms of the $\heta_{ij}$'s, one
obtains $A_\Theta$ as the sum of nine integrals. However, many
of these integrals vanish for purely dimensional reasons. 
After rearranging the points in $C_3(M)$, we can
rewrite $A_\Theta$ as the sum of three contributions:
\[
A_\Theta(M) = A_1(M) + 6\,A_2(M) + 6\,A_3(M),
\]
with
\[
\begin{split}
A_1(M) &= \int_{C_3(M)} \heta_{12}^3\,v_3 = 
\int_{C_2(M)} \heta^3,\\
A_2(M) &= \int_{C_3(M)} \heta_{12}^2\,\heta_{23}\,v_3,\\
A_3(M) &= \int_{C_3(M)} \heta_{12}\,\heta_{23}\,\heta_{31}
\,v_3,
\end{split}
\]
which are graphically represented in fig.\ \ref{figTheta}.
\begin{figure}
\unitlength 1.00mm
\linethickness{0.4pt}
\begin{picture}(111.00,31.00)
\multiput(20.00,30.00)(1.07,-0.12){2}{\line(1,0){1.07}}
\multiput(22.15,29.77)(0.34,-0.12){6}{\line(1,0){0.34}}
\multiput(24.20,29.08)(0.19,-0.11){10}{\line(1,0){0.19}}
\multiput(26.05,27.96)(0.12,-0.11){13}{\line(1,0){0.12}}
\multiput(27.62,26.47)(0.11,-0.16){11}{\line(0,-1){0.16}}
\multiput(28.84,24.68)(0.11,-0.29){7}{\line(0,-1){0.29}}
\multiput(29.64,22.68)(0.12,-0.71){3}{\line(0,-1){0.71}}
\put(29.99,20.54){\line(0,-1){2.16}}
\multiput(29.87,18.38)(-0.12,-0.42){5}{\line(0,-1){0.42}}
\multiput(29.29,16.30)(-0.11,-0.21){9}{\line(0,-1){0.21}}
\multiput(28.28,14.39)(-0.12,-0.14){12}{\line(0,-1){0.14}}
\multiput(26.88,12.74)(-0.16,-0.12){11}{\line(-1,0){0.16}}
\multiput(25.16,11.43)(-0.25,-0.11){8}{\line(-1,0){0.25}}
\multiput(23.19,10.52)(-0.53,-0.12){4}{\line(-1,0){0.53}}
\put(21.08,10.06){\line(-1,0){2.16}}
\multiput(18.92,10.06)(-0.53,0.12){4}{\line(-1,0){0.53}}
\multiput(16.81,10.52)(-0.25,0.11){8}{\line(-1,0){0.25}}
\multiput(14.84,11.43)(-0.16,0.12){11}{\line(-1,0){0.16}}
\multiput(13.12,12.74)(-0.12,0.14){12}{\line(0,1){0.14}}
\multiput(11.72,14.39)(-0.11,0.21){9}{\line(0,1){0.21}}
\multiput(10.71,16.30)(-0.12,0.42){5}{\line(0,1){0.42}}
\put(10.13,18.38){\line(0,1){2.16}}
\multiput(10.01,20.54)(0.12,0.71){3}{\line(0,1){0.71}}
\multiput(10.36,22.68)(0.11,0.29){7}{\line(0,1){0.29}}
\multiput(11.16,24.68)(0.11,0.16){11}{\line(0,1){0.16}}
\multiput(12.38,26.47)(0.12,0.11){13}{\line(1,0){0.12}}
\multiput(13.95,27.96)(0.19,0.11){10}{\line(1,0){0.19}}
\multiput(15.80,29.08)(0.52,0.12){8}{\line(1,0){0.52}}
\put(30.00,20.00){\circle*{2.00}}
\put(10.00,20.00){\circle*{2.00}}
\put(10.00,20.00){\line(1,0){20.00}}
\put(21.00,30.00){\vector(1,0){0.2}}
\put(20.00,30.00){\line(1,0){1.00}}
\put(21.00,10.00){\vector(1,0){0.2}}
\put(20.00,10.00){\line(1,0){1.00}}
\put(21.00,20.00){\vector(1,0){0.2}}
\put(20.00,20.00){\line(1,0){1.00}}
\put(55.00,20.00){\circle{10.00}}
\put(60.00,20.00){\line(1,0){10.00}}
\put(50.00,20.00){\circle*{2.00}}
\put(60.00,20.00){\circle*{2.00}}
\put(70.00,20.00){\circle*{2.00}}
\put(56.00,25.00){\vector(1,0){0.2}}
\put(55.00,25.00){\line(1,0){1.00}}
\put(56.00,15.00){\vector(1,0){0.2}}
\put(55.00,15.00){\line(1,0){1.00}}
\put(66.00,20.00){\vector(1,0){0.2}}
\put(65.00,20.00){\line(1,0){1.00}}
\put(90.00,12.00){\line(1,0){20.00}}
\put(90.00,12.00){\circle*{2.00}}
\put(110.00,12.00){\circle*{2.00}}
\put(99.00,12.00){\vector(-1,0){0.2}}
\put(100.00,12.00){\line(-1,0){1.00}}
\put(70.00,23.00){\makebox(0,0)[cc]{$v$}}
\put(110.00,15.00){\makebox(0,0)[cc]{$v$}}
\multiput(90.00,12.00)(0.12,0.19){84}{\line(0,1){0.19}}
\multiput(100.00,28.00)(0.12,-0.19){84}{\line(0,-1){0.19}}
\put(100.00,28.00){\circle*{2.00}}
\put(96.00,22.00){\vector(1,2){0.2}}
\multiput(95.00,20.00)(0.11,0.22){9}{\line(0,1){0.22}}
\put(106.00,18.00){\vector(1,-2){0.2}}
\multiput(105.00,20.00)(0.11,-0.22){9}{\line(0,-1){0.22}}
\put(8.00,17.00){\makebox(0,0)[cc]{$1$}}
\put(32.00,17.00){\makebox(0,0)[cc]{$2$}}
\put(48.00,17.00){\makebox(0,0)[cc]{$1$}}
\put(62.00,17.00){\makebox(0,0)[cc]{$2$}}
\put(70.00,17.00){\makebox(0,0)[cc]{$3$}}
\put(90.00,9.00){\makebox(0,0)[cc]{$1$}}
\put(110.00,9.00){\makebox(0,0)[cc]{$3$}}
\put(100.00,31.00){\makebox(0,0)[cc]{$2$}}
\put(20.00,5.00){\makebox(0,0)[cc]{$A_1$}}
\put(60.00,5.00){\makebox(0,0)[cc]{$A_2$}}
\put(100.00,5.00){\makebox(0,0)[cc]{$A_3$}}
\end{picture}
\caption{}\label{figTheta}
\end{figure}

\begin{Rem} 
The integral $A_1(M)$ has the same form as
the $\Theta$-invariant in \cite{AS}, \cite{K} and \cite{CT}.
\end{Rem}

\begin{Dig}
The three integrals $A_i$ are not the only possible
combinations containing three $\heta$'s. In fact we can also
consider
\[
A_4(M) = \int_{C_4(M)} v_1\,\heta_{12}\,\heta_{23}\,\heta_{34}
\,v_4.
\]
However, one has the following
\begin{Prop}
For any choice of a metric $g$ and a metric connection $\theta$
involved in the definition of $\heta$,
\[
A_2(M) +2\,A_4(M) =0
\]
if $v=v_g$ is the unit volume form determined by the metric.
\end{Prop}
\begin{proof}
First we notice that, by \eqref{Pgheta}, we can rewrite
\[
A_4(M) = \int_M v_g\,P_g^3v_g.
\]
By \eqref{defPg} we have $dP_gv_g=0$; since $H^2(M)=0$, there exists
a form $\gamma\in\Omega^1(M)$ such that $P_gv_g=d\gamma$. By \eqref{defPg}
we also have
\[
dP_g^2-P_g^2d = P_g\pi_h-\pi_hP_g.
\]
Therefore, we get
\[
A_4(M) = \int_M v_g\,P_g^2d\gamma =
\int_M v_g\,(dP_g^2-P_g\pi_h+\pi_hP_g)\gamma=
\int_M v_g\,P_g\gamma,
\]
since $\pi_hv_g=v_g$ and $\pi_h\gamma=0$.
Notice that this expression is independent of the choice of $\gamma$.
For, if we take $\gamma'\in\Omega^1(M)$ such that $d\gamma'=P_gv$, then
$H^1(M)=0$ implies $\gamma'-\gamma=d\delta$ for some 
$\delta\in\Omega^0(M)$; therefore,
\[
\int_M v_g\,P_g\gamma' - \int_M v_g\,P_g\gamma =
\int_M v_g\,P_gd\delta = -\int_M v_g\,dP_g\delta = 0.
\]

Now we introduce the linear operator $R_g:\Omega^*(M)
\to\Omega^{*+1}(M)$ defined by
\[
R_g = \frac12\,\pi_{2*}(\heta^2\pi_1^*\alpha),
\]
so we can write
\[
A_2(M) = 2\,\int_M R_g P_g v_g.
\]
Following the same lines of the proof of Prop.\ \ref{prop-Pgheta}
and using the second identity of \eqref{bpeta10},
we can show that
\[
dR_g + R_g d = -(\hat v_gP_g+P_g\hat v_g),
\]
where $\hat v_g$ is the operator that acts by multiplication for the
volume form $v_g$.  Therefore,
\[
A_2(M) = 2\,\int_M R_g d\gamma =
-2\,\int_M (dR_g + \hat v_gP_g+P_g\hat v_g)\gamma =
-2\int_M v_g P_g\gamma,
\]
since $\hat v_g\gamma=0$.
\end{proof}
\begin{Rem}
If $\heta$ is so chosen as to represent $\bar P_g=d^*\circ G$,
then both $A_2(M)$ and $A_4(M)$ vanish since $\bar P_g v_g=0$. 
\end{Rem}
End of the digression.
\end{Dig}

\subsection{Proof of Thm.\ \ref{thm-Theta}}\lab{ssec-ptTheta}
We will use here a technique similar to that discussed in
\cite{AS}. That is, we will extend our previous construction
from $C_n(M)$ to $C_n(M)\times I$, where $I$ is a parameter space.

We introduce a parameter $\tau$ ranging over the unit interval $I$ and let
all our quantities---the metric $g$, the metric connection $\theta$
and the unit volume form $v$---depend arbitrarily on $\tau$. Then
$A_\Theta$ will become a function on $I$. More precisely,
we introduce the trivial bundles $C_n(M)\times I$ and denote 
by $\pi$ and $\sigma$ the two projections to $C_n(M)$ and $I$
respectively. Then we define
\[
A_{\Theta,\tau}(M) = \sigma_*(\heta_{123}^3\,v_3),
\]
where now $\heta_{123}v_3$ is seen as a form in $\Omega^9(C_3(M)\times I)$.

As for $v$ we take a representative of the class
in $H^3(M\times I)=H^3(M)$ that satisfies
$\sigma_* v=1$. Notice that, as a form, $v$ belongs to
the completion of
$\Omega^3(M)\otimes\Omega^0(I)\oplus\Omega^2(M)\otimes\Omega^1(I)$.

To let the connection vary on $I$, we consider it as a
connection on the pulled-back
bundle
\[
\pi^{-1}P \xrightarrow{p} M\times I.
\]

Now we will construct $\eta$ as a closed form in 
$\Omega^2(\de C_2(M)\times I)$. As before, we can think of
the sphere bundle $\de C_2(M)\times I\to M\times I$ as
$\pi^{-1}P\times_{SO(3)} S^2$. Consider the commutative
diagram:
\[
\begin{CD}
\pi^{-1}P\times_{SO(3)}S^2 @<{\bar p}<< \pi^{-1}P\times S^2\\
@V{\pi^\de}VV @VV{\bar\pi^\de}V\\
M\times I @<<{p}< \pi^{-1}P
\end{CD}
\]
A form $\bar\eta$ defined as in \eqref{bpeta} will
be a closed, basic form in $\Omega^2(\pi^{-1}P\times S^2)$.
As such, it will be the pullback through $\bar p$ of
a form $\eta\in\Omega^2(\de C_2(M)\times I)$. This form will
satisfy the same properties \eqref{bpeta10} and \eqref{bpetaphi}
described in Prop.\ \ref{thm-eta}. Moreover, we have,
in accordance with Lemma \ref{lem-pie}, the
\begin{Lem}
If $\eta$ is defined as before, then
\[
\pi^\de_*\eta^3 = \frac14\, p_1,
\]
where
\[
p_1 = -\frac1{8\pi^2}\Tr F\wedge F = 
\frac1{4\pi^2}\, F^i F_i
\]
is the first Pontrjagin form on $M\times I$.
\lab{lem-eta3}
\end{Lem}
\begin{proof}
Consider $\bar\eta$ as in \eqref{bpeta}. Since $\omega^2=0$,
\begin{multline*}
(4\pi)^3\bar\pi^\de_*\bar\eta^3 =
\bar\pi^\de_*\{3\,\omega\,[d(\theta^ix_i)]^2 + [d(\theta^ix_i)]^3\}=\\
3\,d\theta^i\,d\theta^j\,\bar\pi^\de_*(\omega\, x_ix_j) -
3\,d\theta^i\,\theta^j\,\theta^k\,\bar\pi^\de_*(x_i\,dx_j\,dx_k).
\end{multline*}
A simple evaluation of these integrals shows that
\[
\bar\pi^\de_*(\omega\, x_ix_j) = \frac43\pi\,\delta_{ij}, \quad\quad
\bar\pi^\de_*(x_i\,dx_j\,dx_k) = \frac43\pi\,\epsilon_{ijk}.
\]
Therefore,
\[
(4\pi)^2p^*\pi^\de_*\eta^3 =
(4\pi)^2\bar\pi^\de_*\bar\eta^3 =
d\theta^i\,d\theta_i-\epsilon_{ijk}d\theta^i\,\theta^j\,\theta^k,
\]
which is equal to $F^iF_i$ by \eqref{defF}. (Notice that
$\Tr\theta^4=0$.)
\end{proof}

The extension of $\eta$ to $\heta\in\Omega^2(C_2(M)\times I)$ and the
definition of $\heta_{123}$ as a representative of $H^2(C_3(M)\times I)
=H^2(C_3(M))$ proceeds as before, by taking an appropriate
$\rho\in\Omega^0(C_2(M)\times I)$ and $\alpha\in\Omega^2(M\times M
\times I)$.
(Notice only that
the involution $T$ and the projections $\pi_i$ and
$\pi_{ij}$ act as the identity on $I$.) In particular,
the properties \eqref{prophetaij}, \eqref{propheta123}
and \eqref{propheta123de} still hold.

Now we are in a position to define $A_{\Theta,\tau}$ and to prove
the following
\begin{Lem}
For an arbitrary dependence of $g$, $\theta$ and $v$ on $I$,
we have
\[
A_{\Theta,1}(M)-A_{\Theta,0}(M) =
\int_I dA_{\Theta,\tau}(M) = -\frac14\, \int_{M\times I} p_1.
\]
\end{Lem}
Then Thm.\ \ref{thm-Theta} follows immediately since
\[
\cs_1(M,f)-\cs_0(M,f) = \int_{M\times I} p_1.
\]

\begin{proof}
We use formula \eqref{dpi*} and get
\[
dA_{\Theta,\tau}(M) = -\sigma_*d(\heta_{123}^3\,v_3) +
\sigma^\de_*\iota_\de^*(\heta_{123}^3\,v_3) =
\int_{\de C_3(M)}\iota_\de^*(\heta_{123}^3\,v_3),
\]
since $\heta_{123}^3\,v_3$ is closed.

We will first consider the principal faces of $\de C_3(M)$, i.e.,
the faces $(12)$, $(23)$ and $(31)$. The last two are immediately
seen to give no contribution since, by (\ref{propheta123de}b)
and (\ref{propheta123de}c), there are no forms depending on the point 1
in the first case and no forms depending
on the point $2$ in the second case.
Therefore, we are left only with the contribution of face $(12)$,
viz.,
\[
-\int_{(12)}\eta_{12}^3 v_3 = -\int_{\de C_2(M)}\eta^3 =
-\frac14\, p_1,
\]
by Lemma \ref{lem-eta3}.

To end our proof we have only to show that the integral
over the hidden face $(123)$ vanishes. This face is a bundle
over $M\times I$,
\[
(123) \xrightarrow{\piott} M\times I,
\]
whose fiber $F$ can be described as
$C_3(\bbR^3)$ modulo translations and scalings. Therefore, $F$
is a 5-dimensional space whose boundary has
three components, denoted by $((12)3)$, $(1(23))$ and
$((31)2)$. A component of $\de F$, say $((12)3)$,
can be then described as follows: fix the translations by $x_3$=0;
so $F$ close to a component
of the boundary looks like $C_2(\bbR^3)$ divided by scalings.
Since $\de C_2(\bbR^3)$ is an $S^2$-bundle over $\bbR^3$, dividing
by the scaling makes each component of $\de F$ an $S^2$-bundle 
over $S^2$.

The integral that we want to evaluate can
be written as
\[
-\int_{(123)} (\eta_{12}+\eta_{23}+\eta_{31})^3\,v_3 =
\int_M v\,\piott_*(\eta_{12}+\eta_{23}+\eta_{31})^3.
\]
We now consider the commutative diagram:
\[
\begin{CD}
\pi^{-1}P\times_{SO(3)}F @<{\bar p}<< \pi^{-1}P\times F\\
@V{\piott}VV @VV{\bar\piott}V\\
M\times I @<<{p}< \pi^{-1}P
\end{CD}
\]
Then, denoting by $\eta_{ij}$, $\omega_{ij}$ and $x_{ij,k}$
the pullbacks of $\eta$, $\omega$ and $x_k$ through the map
$\pi_{ij}$ defined in \eqref{piF},
we have
\begin{multline*}
(4\pi)^3\,p^*\,\piott_*(\eta_{12}+\eta_{23}+\eta_{31})^3 =
(4\pi)^3\,\bar\piott_*\,\bar p^*(\eta_{12}+\eta_{23}+\eta_{31})^3 =\\
\bar\piott_*\left\{
\omega_{12}+\omega_{23}+\omega_{31} +
d\left[\theta^i
\left(x_{12,i}+
x_{23,i}+x_{31,i}
\right)\right]\right\}^3=\\
3\,\theta^i\,\int_F\left(
\omega_{12}+\omega_{23}+\omega_{31} 
\right)^2
\,d\left(
x_{12,i}+x_{23,i}+x_{31,i}
\right)=\\
3\,\theta^i\int_{\de F}
\left(
\omega_{12}+\omega_{23}+\omega_{31} 
\right)^2
\,\left(
x_{12,i}+x_{23,i}+x_{31,i}
\right)=0,
\end{multline*}
The last identity follows from the fact that
\[
\left(
\omega_{12}+\omega_{23}+\omega_{31}
\right)^2
=0
\]
on $\de F$. In fact, on a face, say $((12)3)$, we have
\[
\omega_{12}+\omega_{23}+\omega_{31}
=\omega_{12},
\]
and similarly on the other faces. 
\end{proof}

This concludes the proof of Thm.\ \ref{thm-Theta}.

\subsection{The evaluation of $A_\Theta$ on the 3-sphere}\lab{ssec-S3}
We may think of $S^3$ as the group $SU(2)$. 
Then, given an element $h$,
we take as unit volume form 
\[
v = C\,\Tr(h^{-1}\,dh)^3, \qquad C=\frac1{96\pi^2},
\]
where the trace is taken in the adjoint representation.

In a left and right invariant metric the Levi-Civita
connection is given by
\[
\nabla_{X_i}X_j = \frac12\,[X_i,X_j]
\]
on a left invariant basis of vector fields. This implies that
the connection form on $P$ when pulled back by a left invariant,
orthonormal frame $f_L$ is given by
\[
f_L^*\theta = \frac12\, h^{-1}\,dh.
\]

Consider now the orientation reversing involution
\[
\gamma:\begin{aligned}[t]
S^3 &\to S^3,\\
h &\mapsto h^{-1},
\end{aligned}
\]
and its lifts to $C_2(S^3)$ and $C_3(S^3)$.
With our choice of $v$, we have
\begin{equation}
\gamma^*v=-v.
\lab{g*v}
\end{equation}
Moreover, if we denote by $R$ the adjoint representation which
corresponds to projecting $S^3$ to $SO(3)$, we can write
\begin{equation}
\gamma^*f_L^*\theta=f_L^*\theta^{R(h^{-1})}=
h\,f_L^*\theta\,h^{-1} + h\,dh^{-1}.
\lab{g*theta}
\end{equation}
Let us now consider the action of $\gamma$ on $\de C_2(S^3)$.
On the base we have the action of $\gamma$ defined before; a point
$\bx\in S^2$ is however sent into $-R(h)\,\bx$. In fact, a point
in the tangent bundle is given by $h\,\exp(\bx)$, with $\bx$
in the Lie algebra. Then
\[
\gamma[h\,\exp(\bx)] = \exp(-\bx)\,h^{-1} =
h^{-1}\ h\, \exp(-\bx)\,h^{-1}.
\]
By \eqref{g*theta}, we also have
\[
\gamma^*F = R(h)\, F, \qquad
\gamma^*D\bx = -R(h)\, D\bx.
\]
Therefore, by
\eqref{bpetaTr}, we conclude that
\begin{equation}
\gamma^*\eta=-\eta.
\end{equation}

We can always choose $\rho\in C_2(M)$ to be invariant under
the action of $\gamma$. Then, by \eqref{dreta}, we see that
\[
\gamma^*d\alpha=-d\alpha.
\]
Thus, up to an exact term, we can choose $\alpha$ to be odd and,
finally, obtain
\begin{equation}
\gamma^*\heta = -\heta;
\end{equation}
consequently, we have
\begin{equation}
\gamma^*\heta_{123} = -\heta_{123}.
\end{equation}

We know by Thm.\ \ref{thm-Theta} that the value of $A_\Theta$ does not
depend on these choices, as long as we do not change our connection 
$\theta$. Therefore, with this fixed choice of $\theta$, we have
\begin{equation}
A_\Theta(S^3) = 0
\lab{ATS3}
\end{equation}
since we have found an involution,
$\gamma$, that reverses the orientation of $C_3(M)$ but leaves
$\heta_{123}^3\,v_3$ unchanged.
Therefore, we have
\[
I_\Theta(S^3,f_L) = \frac14\,\cs(S^3,f_L).
\]
Moreover, since
\[
d f_L^*\theta = -2\,f_L^*\theta^2,
\]
we get
\begin{equation}
I_\Theta(S^3,f_L) = \frac1{24\pi^2}\int f_L^*
\Tr\theta^3 = \frac12.
\end{equation}

If we had instead chosen a right invariant section 
$f_R= h^{-1}\,f_L\, h$, then
\[
f^*_R \theta = -\frac12\,dh\,h^{-1},
\]
and we would have obtained
\begin{equation}
I_\Theta(S^3,f_R) = -\frac12.
\end{equation}

\begin{Rem}
The left and right framings are related by the adjoint map
from $S^3$ to $SO(3)$, and hence it has degree 2. So
the corresponding Chern--Simons terms differ by 4
(see Remark \ref{rem-CS}).
As the Chern--Simons terms of these two framings are clearly
opposite in sign, we could have concluded {\em a priori}
that $\cs(S^3,f_{L,R})$ must be $\pm2$.

Note also that the same arguments would have worked for $M=SO(3)$
and would have yielded half the answer for the $\Theta$-invariant.
\end{Rem}

\section{The higher invariants}
\lab{sec-hi}
Our first step is the construction of closed forms
on $C_n(M)$. To do so, we consider the natural projections
\begin{align*}
\pi_i &: C_n(M)\to M, \\
\pi_{ij} &: C_n(M)\to C_2(M), \quad i\not=j,
\end{align*}
and then pull back the volume form $v\in\Omega^3(M)$ 
and the form $\heta\in\Omega^2(C_2(M))$. We will denote
them by $v_i$ and $\heta_{ij}$. The combination $\heta_{ijk}$
defined in \eqref{hetaijk} is now a closed form in $\Omega^2(C_n(M))$
for any triple of distinct indices $ijk$. Of course, not all these forms
are independent. In particular, we notice that
\[
\heta_{ijk} = (-1)^\sigma\,  \heta_{\sigma(jik)},
\]
where $\sigma$ is a permutation. Finally, if $n=2V$,
a product of $3V$ $\heta_{ijk}$'s will be a top form on $C_n(M)$,
while a product of $3V$ $\heta_{ijk}$'s and one volume form will
be a top form on $C_{n+1}(M)$. It is then natural to consider
the relation between these integrals and trivalent graphs. We start with
the following
\begin{Def}[Kontsevich]
In our context the simplest way to describe the graph cohomology
is as follows.
We call a decorated graph a graph with oriented edges and 
numbered vertices (by convention we start the
enumeration by 1). We require edges always to connect distinct
vertices. If two vertices are connected by exactly one edge,
we call that edge regular.
Moreover, denoting by $V$ the number of vertices and
by $E$ the number of edges, we grade the collection of 
decorated graphs by
\begin{equation}
\begin{split}
\ord\Gamma &= E-V,\\
\deg\Gamma &= 2E -3V.
\end{split}
\lab{ord}
\end{equation}
\end{Def}
The $\Theta$-graph has order 1 and degree 0. Examples of decorated
graphs of order 2 are shown in fig.\ \ref{fig-graphs}; $\Gamma_1$ and
$\Gamma_2$ have degree 0, while $\Gamma'$ has degree 1.
\begin{figure}
\unitlength 1.00mm
\linethickness{0.4pt}
\begin{picture}(115.00,40.00)
\multiput(20.00,40.00)(0.99,-0.10){3}{\line(1,0){0.99}}
\multiput(22.97,39.70)(0.36,-0.11){8}{\line(1,0){0.36}}
\multiput(25.83,38.82)(0.22,-0.12){12}{\line(1,0){0.22}}
\multiput(28.45,37.39)(0.13,-0.11){17}{\line(1,0){0.13}}
\multiput(30.74,35.47)(0.12,-0.15){16}{\line(0,-1){0.15}}
\multiput(32.60,33.14)(0.11,-0.22){12}{\line(0,-1){0.22}}
\multiput(33.96,30.48)(0.12,-0.41){7}{\line(0,-1){0.41}}
\multiput(34.77,27.60)(0.11,-1.49){2}{\line(0,-1){1.49}}
\multiput(35.00,24.63)(-0.09,-0.74){4}{\line(0,-1){0.74}}
\multiput(34.62,21.66)(-0.12,-0.35){8}{\line(0,-1){0.35}}
\multiput(33.67,18.83)(-0.11,-0.20){13}{\line(0,-1){0.20}}
\multiput(32.18,16.24)(-0.12,-0.13){17}{\line(0,-1){0.13}}
\multiput(30.20,14.00)(-0.15,-0.11){16}{\line(-1,0){0.15}}
\multiput(27.82,12.20)(-0.24,-0.12){11}{\line(-1,0){0.24}}
\multiput(25.13,10.90)(-0.41,-0.11){7}{\line(-1,0){0.41}}
\multiput(22.24,10.17)(-1.49,-0.07){2}{\line(-1,0){1.49}}
\multiput(19.25,10.02)(-0.74,0.11){4}{\line(-1,0){0.74}}
\multiput(16.30,10.46)(-0.31,0.11){9}{\line(-1,0){0.31}}
\multiput(13.49,11.49)(-0.20,0.12){13}{\line(-1,0){0.20}}
\multiput(10.94,13.04)(-0.13,0.12){17}{\line(-1,0){0.13}}
\multiput(8.75,15.07)(-0.12,0.16){15}{\line(0,1){0.16}}
\multiput(7.01,17.50)(-0.11,0.25){11}{\line(0,1){0.25}}
\multiput(5.78,20.22)(-0.11,0.49){6}{\line(0,1){0.49}}
\put(5.12,23.13){\line(0,1){2.99}}
\multiput(5.04,26.12)(0.10,0.59){5}{\line(0,1){0.59}}
\multiput(5.56,29.06)(0.11,0.28){10}{\line(0,1){0.28}}
\multiput(6.65,31.84)(0.12,0.18){14}{\line(0,1){0.18}}
\multiput(8.27,34.35)(0.12,0.12){18}{\line(0,1){0.12}}
\multiput(10.36,36.49)(0.16,0.11){15}{\line(1,0){0.16}}
\multiput(12.83,38.17)(0.28,0.12){10}{\line(1,0){0.28}}
\multiput(15.58,39.33)(0.74,0.11){6}{\line(1,0){0.74}}
\multiput(60.00,40.00)(0.99,-0.10){3}{\line(1,0){0.99}}
\multiput(62.97,39.70)(0.36,-0.11){8}{\line(1,0){0.36}}
\multiput(65.83,38.82)(0.22,-0.12){12}{\line(1,0){0.22}}
\multiput(68.45,37.39)(0.13,-0.11){17}{\line(1,0){0.13}}
\multiput(70.74,35.47)(0.12,-0.15){16}{\line(0,-1){0.15}}
\multiput(72.60,33.14)(0.11,-0.22){12}{\line(0,-1){0.22}}
\multiput(73.96,30.48)(0.12,-0.41){7}{\line(0,-1){0.41}}
\multiput(74.77,27.60)(0.11,-1.49){2}{\line(0,-1){1.49}}
\multiput(75.00,24.63)(-0.09,-0.74){4}{\line(0,-1){0.74}}
\multiput(74.62,21.66)(-0.12,-0.35){8}{\line(0,-1){0.35}}
\multiput(73.67,18.83)(-0.11,-0.20){13}{\line(0,-1){0.20}}
\multiput(72.18,16.24)(-0.12,-0.13){17}{\line(0,-1){0.13}}
\multiput(70.20,14.00)(-0.15,-0.11){16}{\line(-1,0){0.15}}
\multiput(67.82,12.20)(-0.24,-0.12){11}{\line(-1,0){0.24}}
\multiput(65.13,10.90)(-0.41,-0.11){7}{\line(-1,0){0.41}}
\multiput(62.24,10.17)(-1.49,-0.07){2}{\line(-1,0){1.49}}
\multiput(59.25,10.02)(-0.74,0.11){4}{\line(-1,0){0.74}}
\multiput(56.30,10.46)(-0.31,0.11){9}{\line(-1,0){0.31}}
\multiput(53.49,11.49)(-0.20,0.12){13}{\line(-1,0){0.20}}
\multiput(50.94,13.04)(-0.13,0.12){17}{\line(-1,0){0.13}}
\multiput(48.75,15.07)(-0.12,0.16){15}{\line(0,1){0.16}}
\multiput(47.01,17.50)(-0.11,0.25){11}{\line(0,1){0.25}}
\multiput(45.78,20.22)(-0.11,0.49){6}{\line(0,1){0.49}}
\put(45.12,23.13){\line(0,1){2.99}}
\multiput(45.04,26.12)(0.10,0.59){5}{\line(0,1){0.59}}
\multiput(45.56,29.06)(0.11,0.28){10}{\line(0,1){0.28}}
\multiput(46.65,31.84)(0.12,0.18){14}{\line(0,1){0.18}}
\multiput(48.27,34.35)(0.12,0.12){18}{\line(0,1){0.12}}
\multiput(50.36,36.49)(0.16,0.11){15}{\line(1,0){0.16}}
\multiput(52.83,38.17)(0.28,0.12){10}{\line(1,0){0.28}}
\multiput(55.58,39.33)(0.74,0.11){6}{\line(1,0){0.74}}
\put(9.00,36.00){\line(1,-1){22.00}}
\put(31.00,14.00){\line(0,1){0.00}}
\put(31.00,14.00){\line(0,0){0.00}}
\put(31.00,14.00){\line(0,0){0.00}}
\put(31.00,36.00){\line(-1,-1){22.00}}
\put(20.00,40.00){\vector(1,0){1.00}}
\put(35.00,25.00){\vector(0,-1){1.00}}
\put(20.00,10.00){\vector(-1,0){1.00}}
\put(5.00,25.00){\vector(0,1){1.00}}
\put(25.00,20.00){\vector(1,-1){1.00}}
\put(15.00,20.00){\vector(-1,-1){1.00}}
\put(7.00,38.00){\makebox(0,0)[cc]{1}}
\put(33.00,38.00){\makebox(0,0)[cc]{2}}
\put(33.00,12.00){\makebox(0,0)[cc]{3}}
\put(7.00,12.00){\makebox(0,0)[cc]{4}}
\put(60.00,25.00){\circle{10.00}}
\put(60.00,30.00){\line(0,1){10.00}}
\put(60.00,20.00){\line(0,-1){10.00}}
\put(60.00,35.00){\vector(0,-1){1.00}}
\put(60.00,15.00){\vector(0,-1){1.00}}
\put(55.00,25.00){\vector(0,-1){1.00}}
\put(65.00,25.00){\vector(0,-1){1.00}}
\put(75.00,25.00){\vector(0,-1){1.00}}
\put(45.00,25.00){\vector(0,1){1.00}}
\put(58.00,37.00){\makebox(0,0)[cc]{1}}
\put(62.00,32.00){\makebox(0,0)[cc]{2}}
\put(58.00,18.00){\makebox(0,0)[cc]{3}}
\put(62.00,13.00){\makebox(0,0)[cc]{4}}
\put(20.00,5.00){\makebox(0,0)[cc]{$\Gamma_1$}}
\put(60.00,5.00){\makebox(0,0)[cc]{$\Gamma_2$}}
\multiput(100.00,40.00)(0.99,-0.10){3}{\line(1,0){0.99}}
\multiput(102.97,39.70)(0.36,-0.11){8}{\line(1,0){0.36}}
\multiput(105.83,38.82)(0.22,-0.12){12}{\line(1,0){0.22}}
\multiput(108.45,37.39)(0.13,-0.11){17}{\line(1,0){0.13}}
\multiput(110.74,35.47)(0.12,-0.15){16}{\line(0,-1){0.15}}
\multiput(112.60,33.14)(0.11,-0.22){12}{\line(0,-1){0.22}}
\multiput(113.96,30.48)(0.12,-0.41){7}{\line(0,-1){0.41}}
\multiput(114.77,27.60)(0.11,-1.49){2}{\line(0,-1){1.49}}
\multiput(115.00,24.63)(-0.09,-0.74){4}{\line(0,-1){0.74}}
\multiput(114.62,21.66)(-0.12,-0.35){8}{\line(0,-1){0.35}}
\multiput(113.67,18.83)(-0.11,-0.20){13}{\line(0,-1){0.20}}
\multiput(112.18,16.24)(-0.12,-0.13){17}{\line(0,-1){0.13}}
\multiput(110.20,14.00)(-0.15,-0.11){16}{\line(-1,0){0.15}}
\multiput(107.82,12.20)(-0.24,-0.12){11}{\line(-1,0){0.24}}
\multiput(105.13,10.90)(-0.41,-0.11){7}{\line(-1,0){0.41}}
\multiput(102.24,10.17)(-1.49,-0.07){2}{\line(-1,0){1.49}}
\multiput(99.25,10.02)(-0.74,0.11){4}{\line(-1,0){0.74}}
\multiput(96.30,10.46)(-0.31,0.11){9}{\line(-1,0){0.31}}
\multiput(93.49,11.49)(-0.20,0.12){13}{\line(-1,0){0.20}}
\multiput(90.94,13.04)(-0.13,0.12){17}{\line(-1,0){0.13}}
\multiput(88.75,15.07)(-0.12,0.16){15}{\line(0,1){0.16}}
\multiput(87.01,17.50)(-0.11,0.25){11}{\line(0,1){0.25}}
\multiput(85.78,20.22)(-0.11,0.49){6}{\line(0,1){0.49}}
\put(85.12,23.13){\line(0,1){2.99}}
\multiput(85.04,26.12)(0.10,0.59){5}{\line(0,1){0.59}}
\multiput(85.56,29.06)(0.11,0.28){10}{\line(0,1){0.28}}
\multiput(86.65,31.84)(0.12,0.18){14}{\line(0,1){0.18}}
\multiput(88.27,34.35)(0.12,0.12){18}{\line(0,1){0.12}}
\multiput(90.36,36.49)(0.16,0.11){15}{\line(1,0){0.16}}
\multiput(92.83,38.17)(0.28,0.12){10}{\line(1,0){0.28}}
\multiput(95.58,39.33)(0.74,0.11){6}{\line(1,0){0.74}}
\put(89.00,14.00){\line(2,5){10.40}}
\put(99.40,40.00){\line(2,-5){10.40}}
\put(115.00,25.00){\vector(0,-1){1.00}}
\put(85.00,25.00){\vector(0,1){1.00}}
\put(94.00,26.00){\vector(-1,-1){1.00}}
\put(105.00,26.00){\vector(1,-1){1.00}}
\put(102.00,38.00){\makebox(0,0)[cc]{1}}
\put(92.00,15.00){\makebox(0,0)[cc]{3}}
\put(107.00,15.00){\makebox(0,0)[cc]{2}}
\put(100.00,10.00){\vector(-1,0){1.00}}
\put(100.00,5.00){\makebox(0,0)[cc]{$\Gamma'$}}
\end{picture}
\caption{}\label{fig-graphs}
\end{figure}
Notice that a trivalent graph $\Gamma$ has always degree zero; moreover,
its order can be written as
\[
\ord\Gamma = \frac V2.
\]

To each decorated trivalent graph $\Gamma$ and 3-manifold $M$ we can 
associate a number $A_\Gamma(M)$ given by
\begin{equation}
A_\Gamma(M) \doteq \int_{C_{n+1}(M)} v_0\,\prod_{(ij)\in\Gamma}\heta_{ij0},
\lab{defAGamma}
\end{equation}
where $n=2\ord\Gamma$, and $(ij)$ denotes the oriented edge 
connecting the vertex $i$ to the vertex $j$. Notice that each vertex
carries a number between 1 and $n$, while 0 denotes the $(n+1)$-st
point in $C_{n+1}(M)$.
\begin{Rem}
If we expand the product of $\heta_{ij0}$'s in $A_\Gamma(M)$, the
term coming from choosing $\heta_{ij}$ in each factor will
read
\[
\int_{C_{n+1}(M)} v_0\,\prod_{(ij)\in\Gamma}\heta_{ij} =
\int_{C_n(M)} \prod_{(ij)\in\Gamma}\heta_{ij};
\]
that is, it corresponds to the usual association of a configuration
space integral to a trivalent graph. 
\end{Rem}

\begin{Exa}
Referring to the trivalent decorated graphs of fig.\ \ref{fig-graphs}, 
we have the following integrals:
\[
\begin{split}
A_{\Gamma_1}(M) &= \int_{C_5(M)} v_0\, 
\heta_{120}\,\heta_{230}\,\heta_{340}\,\heta_{410}\,
\heta_{130}\,\heta_{240},\\
A_{\Gamma_2}(M) &= -\int_{C_5(M)} v_0\, 
\heta_{140}^2\,
\heta_{120}\,\heta_{230}^2\,\heta_{340}.
\end{split}
\]
\end{Exa}

In view of the definition of $A_\Gamma(M)$, we give the collection of 
decorated graphs the structure of an algebra over $\bbQ$, 
and extend \eqref{defAGamma} by linearity. We will denote by
$\mathcal{D}$ this algebra.
Moreover, we introduce
the following equivalence relation: if two decorated
graphs $\Gamma$ and $\Gamma'$ have the same underlying 
graph, we set 
\begin{equation}
\Gamma=(-1)^{(p+l)}\,\Gamma'
\lab{eqrel}
\end{equation}
where $p$ is the order of the permutation of the labeling of the vertices
to go from $\Gamma$ to $\Gamma'$, and $l$ is the number of edges
whose orientation must be reversed. Notice that to equivalent
graphs we associate the same number $A_\Gamma(M)$. We will denote
by $\mathcal{D}'$ the algebra of graphs modulo the above equivalence
relation.

Finally, we introduce a coboundary operator $\delta$ on $\mathcal{D}$.
By definition, $\delta\Gamma$ is the 
signed sum of the decorated graphs that are
obtained by contracting a regular edge one at a time in $\Gamma$.
If the new graph is obtained by contracting the oriented edge connecting
the vertex $i$ to the vertex $j$, we relabel the vertices by
letting decrease by one the vertices greater than $\max\{i,j\}$
and denote by $\min\{i,j\}$ the vertex where the contraction has happened.
Moreover, associate to this contraction a sign $\sigma(i,j)$
defined by
\begin{equation}
\sigma(i,j) = \begin{cases}
(-1)^j &\text{if $j>i$,}\\
(-1)^{i+1} &\text{if $j<i$}.
\end{cases}
\lab{sigmaij}
\end{equation}

\begin{Prop}
The operator $\delta$ descends to $\mathcal{D}'$ and satisfies
$\delta^2=0$ there. Moreover, if we denote by $\mathcal{D}'_{n,t}$
the (equivalence classes of) decorated graphs of order $n$ and
degree $t$, we have
\[
\delta : \mathcal{D}'_{n,t} \to \mathcal{D}'_{n,t+1}.
\]
\lab{prop-delta}
\end{Prop}
\begin{proof}
If we exchange $i$ and $j$ or reverse the arrow connecting them,
we get a minus sign from \eqref{eqrel}. In both cases the roles
of $i$ and $j$ are exchanged. However, $\sigma(i,j)=-\sigma(j,i)$.
Therefore, it does not matter if we exchange $i$ and $j$ first and 
then apply $\delta$ to the edge $(ij)$ or vice versa.

Then consider three vertices $i$, $j$ and $k$. We want to prove
that the exchange of $j$ with $k$ does not interfere with the action
of $\delta$ on $(ij)$. By the previous step, we can assume $i<j$
and $(ij)$ oriented from $i$ to $j$. First suppose $k>i$. We can
also assume $k>j$.
If we contract $(ij)$ we get a factor $(-1)^j$. If we exchange
$j$ and $k$ we get a factor $-1$ and then we have to contract $(ik)$
with a factor $(-1)^k$. Now we have to consider what happens to the 
relabeling. In the first case all indices 
greater than $j$ are reduced 
by one, in the second only those greater than $k$. 
The vertices labeled as $j(j+1)\dots(k-2)(k-1)$ in the first case
are labeled as $(j+1)(j+2)\dots(k-1)j$ in the second. 
Since the two strings have length $k-j$ and are related by one
cyclic rotation, we get a factor $(-1)^{k-j+1}$ if we want to
turn one graph into the other. In summary, if we contract
$(ij)$ we get a sign $(-1)^j$, while if we echange $j$ and $k$
and then contract $(ik)$ we get a sign $(-1)^{k+1}$ and a labeling
that is related to the previous one by the sign $(-1)^{k-j+1}$.
Similarly, we can treat the case $k<i$. Therefore, $\delta$ descends
to $\mathcal{D}'$.

To prove that $\delta^2=0$ on $\mathcal{D}'$,
we check that contracting two different
pairs $(ij)$ and $(rs)$ in the opposite order gives opposite signs.
First assume $j\not=s$ and $i\not=r$. 
By reordering the vertices, we can assume
that $i<j$, $r<s$ and $j<s$. Then contracting $i$ with $j$ gives
$(-1)^j$ and let $s$ (and possibly $r$) decrease by one, so 
contracting $r$ with $s$ gives $(-1)^{s-1}$; if we instead contract
$r$ with $s$ first and then $i$ with $j$, we get  $(-1)^s(-1)^j$
since $j$ is not reduced by one.
If $s=j$ and $i\not=r$, then the two orders in which we can 
contract the pairs $(ij)$ 
and $(rj)$ are related by exchanging $i$ with $r$, which changes
sign on $\mathcal{D}'$. Similarly we treat the case $s\not=j$,
$i=r$.

Finally, we observe that a contraction decreases by one both
the number of vertices $V$ and the number of edges $E$ (remember
that we contract only regular edges). The last claim is
thus a consequence of \eqref{ord}.
\end{proof}

\begin{Exa}
The $\Theta$-graph is a cocycle since it has no regular edges.
\end{Exa}

\begin{Exa}
Referring to the decorated graphs in fig.\ \ref{fig-graphs}, we have
\[
\begin{split}
\delta\Gamma_1 &= 6\,\Gamma',\\
\delta\Gamma_2 &= 2\,\Gamma'.
\end{split}
\]
Therefore, the combination
\[
\Gamma=
-\frac1{12}\,\Gamma_1 + \frac14\,\Gamma_2
\]
is a cocycle. \lab{coGamma}
\end{Exa}

Notice that the action of $\delta$ can be restricted
to the algebra of connected graphs. 
Then we have the following 
\begin{Thm}
If\/ $\Gamma$ is a connected, trivalent cocycle in graph cohomology, 
then there exists
a constant $\phi(\Gamma)$ such that
\[
I_\Gamma(M,f) = A_\Gamma(M) +\phi(\Gamma)\,\cs(M,f)
\]
is an invariant for the framed rational homology 3-sphere $M$.
Moreover, if\/ $\ord\Gamma$ is even, then $\phi(\Gamma)=0$.
\lab{thm-Gamma}
\end{Thm}

{}From Thms. \ref{thm-Theta}, \ref{thm-Gamma} and the discussion
in subsection \ref{ssec-S3}, we get the following
\begin{Cor}
If\/ $\Gamma$ is a connected, trivalent cocycle in graph cohomology, then
the quantity
\[
J_\Gamma(M) = A_\Gamma(M) - 4\,\phi(\Gamma)\, A_\Theta(M)
\]
is an invariant for the rational homology 3-sphere $M$. 
Moreover, if\/ 
$\ord\Gamma$ is even, then $\phi(\Gamma)=0$; if\/ $\ord\Gamma$
is odd, then $J_\Gamma(S^3)=0$.
\end{Cor}
The last statement follows from the fact that, if we choose
the connection $\theta$ as in subsection \ref{ssec-S3}, then
$A_\Theta(S^3)=0$; moreover, the involution $\gamma$ we defined
there reverses the orientation of $C_{n+1}(S^3)$ since $n$ is even,
and sends each $\heta_{ij0}$ and the unit volume form 
into minus themselves.
Since the number of forms $\heta_{ij0}$ is equal to the number of
edges of $\Gamma$, that is, to $3\ord\Gamma$, we see that the integrand
is even when $\ord\Gamma$ is odd.

\subsection{Proof of Thm.\ \ref{thm-Gamma}}\lab{ssec-ptGamma}
As in the proof of Thm.\ \ref{thm-Theta} we introduce the
unit interval $I$ on which all our quantities depend. 
In the following we will use notation similar to that of subsection 
\ref{ssec-ptTheta}. In particular, we still denote by
$\pi$ the projection $M\times I\to M$ and by $\pi^{-1}P$ the pulled-back
orthonormal frame bundle.

If $\mathcal{S}$ is a connected component of codimension one 
of $\de C_n(M)\times I$ describing the coincidence of $p$ points,
we will denote by $\pi^\de$ its projection to $C_{n-p+1}(M)\times I$
and by $F$ its fiber. This fiber 
is given by $C_p(\bbR^3)$ modulo translations and rescalings.
Therefore, it is a $(3p-4)$-manifold with corners.
Then we have the following
\begin{Lem}[Kontsevich \cite{K}]
If $F$ is the fiber of the face $\mathcal{S}$
associated with the collision of $q$ points $\bx_1,\dots,\bx_q$,
and by $\omega_{ij}$
we denote the pullbacks of the volume form of the sphere
through the projections $\pi_{ij}$ defined in \eqref{piF}, 
then, for any triple of indices $i,j,k$ ($i\not=j,i\not=k$),
\[
\int_{\bx_i} \omega_{ij}\,\omega_{ik} = 0.
\]
\lab{lem-omega}
\end{Lem}
\begin{proof}
If $j=k$ the identity follows since $\omega_{ij}^2=0$.

If $j\not=k$, then we have at least three vertices in $F$.
Then the integration on $\bx_i$ extends to all $\bbR^3$
with some
points blown up since we can use the other two vertices
to fix translations and scalings.
The form $\omega_{ij}\,\omega_{ik}$
is regular except at the points $\bx_i=\bx_j$ and 
$\bx_i=\bx_k$; so
we can extend the integration of $\bx_i$ to $\bbR^3$ with only these two 
points blown up. Then we consider the orientation reversing involution
\[
\bx_i'=\bx_j+\bx_k-\bx_i.
\]
The forms $\omega_{ij}$ and $\omega_{ik}$ depend only
on the difference of $\bx_i$ with $\bx_j$ and $\bx_k$, so
they are sent to $\omega_{ki'}=-\omega_{i'k}$ and 
$\omega_{ji'}=-\omega_{i'j}$ respectively. 
Since this involution reverses the orientation of the manifold
and leaves the integrand form unchanged, the integral must vanish.
\end{proof}

Let us denote by $\pi_x$ the projection $C_{n-p+1}(M)\to M$ 
corresponding to the point where all the $p$ points collapsed.
Then we define $\tilde\pi=\pi\circ\pi_x$ and consider the pulled-back
bundle $\tilde\pi^{-1}P=\pi_x^{-1}\pi^{-1}P$. 
We can now identify
$\mathcal{S}$ with $\tilde\pi^{-1}P\times_{SO(3)}F$ and consider
the commutative diagram:
\[
\begin{CD}
\tilde\pi^{-1}P\times_{SO(3)}F @<{\bar p}<< \tilde\pi^{-1}P\times F\\
@V{\pi^\de}VV @VV{\bar\pi^\de}V\\
C_{n-p+1}(M)\times I @<<{p}< \tilde\pi^{-1}P
\end{CD}
\]
The form $\eta_{ij}$ on this boundary is the
pullback of the form $\eta$ on $\de C_2(M)\times I=
\pi^{-1}P\times_{SO(3)}S^2$
through the projection $\pi_{x,ij}$ given by the composition of
the projection $\pi_{ij}: F\to S^2$ defined in \eqref{piF} with
the projection $\pi_x$. With this notation we can state the following
\begin{Lem}[Axelrod and Singer \cite{AS}]
If $\lambda$ is a form in 
$\Omega^*(\mathcal{S})$ given by a sum of products of forms
$\eta_{ij}$, then
$\pi_*^\de\lambda$ is the pullback through $\pi_x$ of
a characteristic form on $M\times I$.
In particular, it vanishes unless it is of degree zero, in which
case it is a constant, or of degree four, in which case it is a
multiple of the first Pontrjagin form $p_1$.
\lab{lem-p1}
\end{Lem}
\begin{proof}
By construction $\bar p^*\lambda$ is a polynomial in $\pi_x^*\theta$ and
$\pi_x^*d\theta$ with coefficients in $\Omega^*(F)$. Therefore,
$p^*\pi^\de_*\lambda=\bar\pi^\de_*\bar p^*\lambda$ is a polynomial
in  $\pi_x^*\theta$ and $\pi_x^*d\theta$ and is basic. This means
that it is a characteristic form on $C_{n-p+1}(M)\times I$ obtained
by pullback from $M\times I$ through $\pi_x$. Since $M\times I$ is
4-dimensional, its only characteristic forms are the constant function
and the first Pontrjagin form.
\end{proof}

Now we introduce the projection $\sigma$ from $C_{n+1}(M)\times I$
to $I$, and, for each decorated trivalent graph, we define
\[
A_{\Gamma,\tau}(M) = \sigma_*\left(
v_0\,\prod_{(ij)\in\Gamma}\heta_{ij0}
\right).
\]
Since the integrand is a closed form, $dA_{\Gamma,\tau}$ will
be given only by boundary contributions. We have the following
\begin{Lem}
The only boundary terms that contribute to $dA_{\Gamma,\tau}$
come from:
\begin{enumerate}
\item faces corresponding to the collapse of two points connected
by exactly one edge in $\Gamma$;
\item the face corresponding to the collapse of
all points but the point labeled 0; in this case, the contribution
is a multiple of the first Pontrjagin form if\/ $\ord\Gamma$ is odd,
and vanishes if\/ $\ord\Gamma$ is even.
\end{enumerate}\lab{lem-Gamma}
\end{Lem}

Assuming this lemma for the moment,
suppose that $1$ and $2$ are points as in case 1. Then, on
the face where $1$ and $2$ come together, $\heta_{120}=-\eta_{12}$
and its pushforward gives $-1$. 
Notice that the rest of the graph
is left unchanged.
Now compare what happens if we contract the pair
$(ij)$, assuming $i<j$, instead of $(12)$. 
To do so we bring the pair $i$ to $1$ and $j$ to $2$ by using
cyclic rotations, then contract and
rotate again to bring the vertex $i$ into its original position.
More explicitly, first we rotate
the points $12\dots i\dots j\dots n$ to 
$i(i+1)\dots j\dots n12\dots(i-1)$. 
Since we have done $i-1$ cyclic rotations in an even chain, we get
the sign $(-1)^{i-1}$. Then we rotate $(i+1)\dots j$ to
$j(i+1)\dots(j-1)$. This gives sign $(-1)^{j-i+1}$ since we have
done one cyclic rotation in a chain of $j-i$ elements. Then
we contract the two vertices.
Finally, we rotate $i$ back to its original position; this gives no sign
since the chain is odd now. Thus, the contraction of $i$ with $j$,
with $i<j$, has sign $(-1)^j$. This is in accordance with the
sign convention \eqref{sigmaij} in the definition of the coboundary
operator $\delta$.
Therefore, we have proved the following
\begin{Cor}
If\/ $\Gamma$ is a connected, trivalent cocycle in graph cohomology,
for an arbitrary dependence of $g$, $\theta$ and $v$ on $I$,
\[
A_{\Gamma,1}(M)-A_{\Gamma,0}(M) =
\int_I dA_{\Gamma,\tau}(M) = -\phi(\Gamma) \int_{M\times I} p_1,
\]
where $\phi(\Gamma)$ is a number that depends only on the
cocycle\/ $\Gamma$ and vanishes if\/ $\ord\Gamma$ is even.
More generally, if\/ $\Gamma$ is not a cocycle, we have
\[
dA_{\Gamma,\tau}(M) = -A_{\delta\Gamma,\tau}(M)
-\phi(\Gamma) \int_{M\times I} p_1.
\]
\end{Cor}
{}From this corollary we obtain Thm.\ \ref{thm-Gamma}.

\begin{proof}[Proof of Lemma \ref{lem-Gamma}]
To fix our notation, in the following we will denote by $\mathcal{S}$
a boundary face and by $\pi^\de$ its projection to $C_{n-q+2}(M)\times I$,
where $q$ ($q>1$)
is the number of collapsing points. We will denote by $F$
the $(3q-4)$-dimensional fiber. We will denote by $\lambda$ the
restriction to the boundary of one summand of
the integrand form in $A_{\Gamma,\tau}$,
and will write
\[
\lambda = \lambda_1\,\pi^{\de*}\lambda_2.
\]

First we show that the boundary faces involving the point labeled by 0
do not contribute. 
Suppose that, besides 0, we have $p$ other points ($p>0$)
coming together (that is, we are taking $q=p+1$). 
The boundary then projects to $C_{n-p+1}$ and the fiber
$F$ has dimension $3p-1$. Denote by $e$ the number of edges in $\Gamma$
connecting two points on the boundary face
$\mathcal{S}$ and by $e_0$ the number of edges in $\Gamma$ connecting
a point in $\mathcal{S}$ to a point outside $\mathcal{S}$. Since
$\Gamma$ is a trivalent graph, we have
\[
2e+e_0=3p.
\]
Notice that, if both $i$ and $j$ belong to $\mathcal{S}$, then
$\heta_{ij0}=-(\eta_{ij}+\eta_{j0}+\eta_{0i})$; if only $i$
belongs to $\mathcal{S}$, then $\heta_{ij0}=-\eta_{i0}$. 
Therefore, each
edge with at least one vertex in $\mathcal{S}$ contributes
with a 2-form in the vertical direction; that is,
\[
\deg\lambda_1 = 2(e+e_0) = 3p+e_0.
\]
Since the fiber has dimension $3p-1$, we see that
\[
\deg\pi^\de_*\lambda_1 = e_0+1.
\]
Thus, if $e_0>0$, $v_0\,\pi^\de_*\lambda_1=0$ since $M\times I$ has
dimension four. On the other hand, if $e_0=0$, then
$\pi^\de_*\lambda_1$ vanishes by Lemma \ref{lem-p1}.

Next we come to the case where $q$ points ($1<q<n+1$) come together
and the point labeled by 0 is not involved. In this case, if both $i$
and $j$ belong to $\mathcal{S}$, then $\heta_{ij0}=-\eta_{ij}$. If, however,
at least one vertex does not belong to $\mathcal{S}$, then $\heta_{ij0}$
will be basic (that is, it will contribute to $\lambda_2$).
Therefore, the form degree in the vertical direction is equal
to the number $e$ of edges connecting vertices inside $\mathcal{S}$.
Again we have the relation
\[
2e+e_0 = 3q,
\]
which implies
\[
\deg\lambda_1 = 3q-e_0.
\]
Since the fiber has dimension $3q-4$, we get
\[
\deg\pi^\de_*\lambda_1 = 4-e_0.
\]
By Lemma \ref{lem-p1}, we see that we have a nonvanishing contribution
only if $e_0=4$ or $e_0=0$; moreover, in the latter case,
$\pi^\de_*\lambda_1$ is a multiple of the first Pontrjagin form.
Notice that, since we are considering only connected diagrams,
this case corresponds to case 2 in the Lemma we are proving.
To prove that this contribution vanishes if $\ord\Gamma$ is even,
consider the involution that reverses all coordinates $\bx_i$ in
$F$. Since $n$ is even, this involution is orientation reversing.
(Notice that we can represent $F$ as $S^{3n-4}$, with some
submanifolds blown up, and that the involution corresponds
to the antipodal map.)
On the other hand, each $\eta_{ij}$ is sent into $-\eta_{ij}$, and
since the number of $\eta$'s is $E=3\ord\Gamma$, we see that the 
integrand
does not change sign if $\ord\Gamma$ is even. In this case, then,
the integral vanishes.

We are now left with case $e_0=4$. Notice that, since 
$\deg\lambda_1=\dim F$ now, we must select the top form
on $F$ in $\lambda_1$; that is, we must
replace each $\eta_{ij}$ with $\omega_{ij}/(4\pi)$.
We have two possibilities:
\begin{enumerate}
\item there is at least one vertex in $\mathcal{S}$ connected
to a vertex outside $\mathcal{S}$ through exactly one leg;
\item there are two vertices in $\mathcal{S}$ each of which
is connected with the outside through exactly two legs (the two
legs can connect the vertex inside with the same vertex outside or 
with two distinct vertices outside), and no other
vertex in $\mathcal{S}$ is connected to the outside. 
\end{enumerate}
In case 1, we can apply Lemma \ref{lem-omega} and conclude that
$\pi^\de_*\lambda_1=0$. In case 2, we notice that the two
vertices under consideration can be connected by an edge
only if $q=2$ since the diagram is connected, and this corresponds
to case 1 in the Lemma we are proving. Thus, we assume that the
two vertices, which we denote by $i$ and $j$,
are not connected and $q>2$. In this case, there exists
a third vertex $k$ that is conncted to $i$ through exactly one
edge to which we associate the form $\omega_{ik}$.
Then we integrate $\omega_{ik}$ over the position of the point
$i$. If we make the change of variables
\[
\bx_i' = \bx_i-\bx_k,
\]
we realize that the result of this integration does not depend
on $\bx_k$. Thus, we can see the vertex $k$ as if it were not
connected to $i$. It will however be connected to two other vertices,
possibly not distinct, in $\mathcal{S}$. Then we can use
Lemma \ref{lem-omega} and conclude that
$\pi^\de_*\lambda_1$ vanishes also in this case.
\end{proof}

This concludes the proof of Thm.\ \ref{thm-Gamma}.

\section{Knots in a rational homology 3-sphere}
The forms $\heta_{ijk}$ we have introduced in the previous sections
allow for the construction of invariants of knots 
$K\subset M$---when
$M$ is a rational homology 3-sphere---generalizing the case of knots
in $\bbR^3$ discussed in \cite{BT}.

In general (for details see \cite{BT}), an imbedding
\[
f:X\hookrightarrow Y
\]
gives rise to natural imbeddings
\begin{align*}
f^n :  X^n &\hookrightarrow Y^n,\\
\intertext{and}
C_n^f :  C_n(X) &\hookrightarrow C_n(Y).
\end{align*}
Moreover, since we have natural projections 
$\pi:C_{n+t}(Y)\to C_n(Y)$, we can consider the pulled-back
bundles 
\[
\mathcal{C}_{n,t}^f=(C_n^f)^{-1}C_{n+t}(Y)\to C_n(X).
\]
We will then have natural projections 
$\mathcal{C}_{n,t}^f\to \mathcal{C}_{n-r,t-s}^f$.

In our case we set $X=S^1$ and $Y=M$, and consider a family
of imbeddings
\[
K_\tau:S^1\hookrightarrow M, \qquad \tau\in I,
\]
where $I$ is the unit interval. 
{}From this we can define the families of imbeddings
$C_{n}^{K_\tau}$ and of
bundles $\mathcal{C}_{n,t}^{K_\tau}$. 

Notice however that
$C_n(S^1)$ has $n!$ distinct connected components. We pick up
one of them by choosing a fixed
ordering of the points on $S^1$. This connected
configuration space will be denoted by $\wC_n(S^1)$. Correspondingly
we will have the families of imbeddings $\wC_{n}^{K_\tau}$ 
and of bundles $\wcC_{n,t}^{K_\tau}$.

Then we recall that we are also interested in varying the metric, 
the connection form and the unit volume form on $M$. 
We can take the parameter
$\tau$ to belong to the same unit interval $I$,
so we are led to consider the map
\[
\hat K : \begin{array}[t]{ccc}
S^1\times I &\to& M\times I,\\
(\alpha,\tau) &\mapsto& (K_\tau(\alpha),\tau),
\end{array}
\]
and its generalizations
\[
\hwC_n^K : \wC_n(S^1)\times I \to C_n(M)\times I. 
\]
Finally, the natural projections $\pi:C_{n+t}(M)\times I\to C_n(M)\times I$
allow us to define the pulled-back bundles
\[
\hwcC_{n,t}^K = (\hwC_n^K)^{-1}(C_{n+t}(M)\times I).
\]
Again we have natural projections $\hwcC_{n,t}^K\to\hwcC_{n-r,t-s}^K$;
the case $r=n, t=s$ yields the projection
\[
\sigma: \hwcC_{n,t}^K\to I.
\]

Using the maps
$\wC_n^K$ we can pull back the forms 
$\heta_{ij}\in\Omega^2(C_n(M))$. We will
keep denoting them by $\heta_{ij}$ to avoid cumbersome notation.
Similarly we can pull back the forms 
$\heta_{ij}\in\Omega^2(C_n(M)\times I)$ 
by the maps $\hwC_n^K$.

Notice
that a form $\heta_{ij}\in\Omega^2(\hwcC_{n,t}^K)$ 
(more precisely we should write
$\hwC_n^{K*}\heta_{ij}$) depends on $\tau$ in two ways:
through the metric connection $\theta$ and through the map $\hwC_n^K$.

As in the case of manifold invariants, we will look for configuration
space integrals that yield functions on $I$. A constant function
will then be a knot invariant.

The simplest quantity we can write down with these ingredients
is the {\em self-linking number}
\begin{equation}
\sln(K,M) \doteq \int_{\wcC_{2,0}^K} \heta,
\lab{sln}
\end{equation}
which is {\em not} a knot invariant. In fact, 
given a family of imbeddings $K_\tau$, we can write
\[
\sln(K_\tau,M) = \sigma_*\heta,
\]
where $\sigma$ is the projection $\hwcC_{2,0}^K\to I$.
We then note that 
$\hwC_2^{K*}v_i=0$, so $\heta$ is a closed form on
$\hwcC_{2,0}^K$; moreover,
$\de \wcC_{2,0}^K = S^1\times\{-1,1\}$. Therefore, we have
\begin{equation}
d\sln(K_\tau,M) = \sigma_*^\de\heta =
\int_{\de \wcC_{2,0}^{K_\tau}} \heta = 
-2\int_{S^1} \psi_{K_\tau}^*\eta.
\lab{dsln}
\end{equation}
For a given imbedding $K$, the map 
\begin{equation}
\psi_K : S^1 \to \de C_2(M)
\lab{defpsi}
\end{equation}
is defined as follows: Consider the tangent map $K_*:TS^1\to TM$ and
its restriction to the sphere bundles. Since $S(TS^1)=S^1
\times\{-1,1\}$,
we actually have two maps---opposite to each other---from 
$S^1$ to $S(TM)=\de C_2(M)$, one
corresponding to the point $-1$ and the other to the point $1$.
We take $\psi_K$ as the latter.

In the next subsection we will show how to associate knot invariants
to cocycles in a new graph cohomology. We will see that the only
possible failure for the integrals we write down 
to be true invariants is given by a term
proportional to $d\sln(K_\tau,M)$. Therefore, subtracting the correct
multiple of the self-linking number, we will get knot invariants.

\subsection{Knot invariants}
To define knot invariants, we have to introduce an appropriate graph
cohomology. 
Essentially we use the same diagrams considered in Sec.\ \ref{sec-hi},
but with some important modifications.

\begin{Def}
We call a decorated graph for knots a decorated graph with
a distinguished loop (representing the knot)
on which we orient all the edges consistently. We call {\em external}
the vertices on the knot and {\em internal} the others. We assume
that there are always at least two external vertices.
We call internal the edges which do not constitute the knot and
external those which do. Following \cite{BT}, we call 
{\em connected} a decorated graph for knots
such that its underlying graph is connected after removing
any pair of external edges. (In \cite{AF}, such a graph
is called prime.)
Finally, denoting by $E$ the number
of internal edges and by $V_i$ and $V_e$ the number of internal
and external vertices, we grade the collection of decorated
graphs for knots by
\begin{equation}
\begin{split}
\ord\Gamma &= E-V_i,\\
\deg\Gamma &= 2E-3V_i-V_e.
\end{split}
\end{equation}
\end{Def}

The $\Theta$-graph can be seen as a connected decorated graph for knots
of order 1 and degree 0
if its outer circle is reinterpreted as the knot.
Examples of connected decorated graphs for knots of order
2 and degree 0 are shown in fig.\ \ref{fig-graphsK};
\begin{figure}
\unitlength 1.00mm
\linethickness{0.4pt}
\begin{picture}(115.00,40.00)
\multiput(20.00,40.00)(0.99,-0.10){3}{\line(1,0){0.99}}
\multiput(22.97,39.70)(0.36,-0.11){8}{\line(1,0){0.36}}
\multiput(25.83,38.82)(0.22,-0.12){12}{\line(1,0){0.22}}
\multiput(28.45,37.39)(0.13,-0.11){17}{\line(1,0){0.13}}
\multiput(30.74,35.47)(0.12,-0.15){16}{\line(0,-1){0.15}}
\multiput(32.60,33.14)(0.11,-0.22){12}{\line(0,-1){0.22}}
\multiput(33.96,30.48)(0.12,-0.41){7}{\line(0,-1){0.41}}
\multiput(34.77,27.60)(0.11,-1.49){2}{\line(0,-1){1.49}}
\multiput(35.00,24.63)(-0.09,-0.74){4}{\line(0,-1){0.74}}
\multiput(34.62,21.66)(-0.12,-0.35){8}{\line(0,-1){0.35}}
\multiput(33.67,18.83)(-0.11,-0.20){13}{\line(0,-1){0.20}}
\multiput(32.18,16.24)(-0.12,-0.13){17}{\line(0,-1){0.13}}
\multiput(30.20,14.00)(-0.15,-0.11){16}{\line(-1,0){0.15}}
\multiput(27.82,12.20)(-0.24,-0.12){11}{\line(-1,0){0.24}}
\multiput(25.13,10.90)(-0.41,-0.11){7}{\line(-1,0){0.41}}
\multiput(22.24,10.17)(-1.49,-0.07){2}{\line(-1,0){1.49}}
\multiput(19.25,10.02)(-0.74,0.11){4}{\line(-1,0){0.74}}
\multiput(16.30,10.46)(-0.31,0.11){9}{\line(-1,0){0.31}}
\multiput(13.49,11.49)(-0.20,0.12){13}{\line(-1,0){0.20}}
\multiput(10.94,13.04)(-0.13,0.12){17}{\line(-1,0){0.13}}
\multiput(8.75,15.07)(-0.12,0.16){15}{\line(0,1){0.16}}
\multiput(7.01,17.50)(-0.11,0.25){11}{\line(0,1){0.25}}
\multiput(5.78,20.22)(-0.11,0.49){6}{\line(0,1){0.49}}
\put(5.12,23.13){\line(0,1){2.99}}
\multiput(5.04,26.12)(0.10,0.59){5}{\line(0,1){0.59}}
\multiput(5.56,29.06)(0.11,0.28){10}{\line(0,1){0.28}}
\multiput(6.65,31.84)(0.12,0.18){14}{\line(0,1){0.18}}
\multiput(8.27,34.35)(0.12,0.12){18}{\line(0,1){0.12}}
\multiput(10.36,36.49)(0.16,0.11){15}{\line(1,0){0.16}}
\multiput(12.83,38.17)(0.28,0.12){10}{\line(1,0){0.28}}
\multiput(15.58,39.33)(0.74,0.11){6}{\line(1,0){0.74}}
\multiput(60.00,40.00)(0.99,-0.10){3}{\line(1,0){0.99}}
\multiput(62.97,39.70)(0.36,-0.11){8}{\line(1,0){0.36}}
\multiput(65.83,38.82)(0.22,-0.12){12}{\line(1,0){0.22}}
\multiput(68.45,37.39)(0.13,-0.11){17}{\line(1,0){0.13}}
\multiput(70.74,35.47)(0.12,-0.15){16}{\line(0,-1){0.15}}
\multiput(72.60,33.14)(0.11,-0.22){12}{\line(0,-1){0.22}}
\multiput(73.96,30.48)(0.12,-0.41){7}{\line(0,-1){0.41}}
\multiput(74.77,27.60)(0.11,-1.49){2}{\line(0,-1){1.49}}
\multiput(75.00,24.63)(-0.09,-0.74){4}{\line(0,-1){0.74}}
\multiput(74.62,21.66)(-0.12,-0.35){8}{\line(0,-1){0.35}}
\multiput(73.67,18.83)(-0.11,-0.20){13}{\line(0,-1){0.20}}
\multiput(72.18,16.24)(-0.12,-0.13){17}{\line(0,-1){0.13}}
\multiput(70.20,14.00)(-0.15,-0.11){16}{\line(-1,0){0.15}}
\multiput(67.82,12.20)(-0.24,-0.12){11}{\line(-1,0){0.24}}
\multiput(65.13,10.90)(-0.41,-0.11){7}{\line(-1,0){0.41}}
\multiput(62.24,10.17)(-1.49,-0.07){2}{\line(-1,0){1.49}}
\multiput(59.25,10.02)(-0.74,0.11){4}{\line(-1,0){0.74}}
\multiput(56.30,10.46)(-0.31,0.11){9}{\line(-1,0){0.31}}
\multiput(53.49,11.49)(-0.20,0.12){13}{\line(-1,0){0.20}}
\multiput(50.94,13.04)(-0.13,0.12){17}{\line(-1,0){0.13}}
\multiput(48.75,15.07)(-0.12,0.16){15}{\line(0,1){0.16}}
\multiput(47.01,17.50)(-0.11,0.25){11}{\line(0,1){0.25}}
\multiput(45.78,20.22)(-0.11,0.49){6}{\line(0,1){0.49}}
\put(45.12,23.13){\line(0,1){2.99}}
\multiput(45.04,26.12)(0.10,0.59){5}{\line(0,1){0.59}}
\multiput(45.56,29.06)(0.11,0.28){10}{\line(0,1){0.28}}
\multiput(46.65,31.84)(0.12,0.18){14}{\line(0,1){0.18}}
\multiput(48.27,34.35)(0.12,0.12){18}{\line(0,1){0.12}}
\multiput(50.36,36.49)(0.16,0.11){15}{\line(1,0){0.16}}
\multiput(52.83,38.17)(0.28,0.12){10}{\line(1,0){0.28}}
\multiput(55.58,39.33)(0.74,0.11){6}{\line(1,0){0.74}}
\multiput(100.00,40.00)(0.99,-0.10){3}{\line(1,0){0.99}}
\multiput(102.97,39.70)(0.36,-0.11){8}{\line(1,0){0.36}}
\multiput(105.83,38.82)(0.22,-0.12){12}{\line(1,0){0.22}}
\multiput(108.45,37.39)(0.13,-0.11){17}{\line(1,0){0.13}}
\multiput(110.74,35.47)(0.12,-0.15){16}{\line(0,-1){0.15}}
\multiput(112.60,33.14)(0.11,-0.22){12}{\line(0,-1){0.22}}
\multiput(113.96,30.48)(0.12,-0.41){7}{\line(0,-1){0.41}}
\multiput(114.77,27.60)(0.11,-1.49){2}{\line(0,-1){1.49}}
\multiput(115.00,24.63)(-0.09,-0.74){4}{\line(0,-1){0.74}}
\multiput(114.62,21.66)(-0.12,-0.35){8}{\line(0,-1){0.35}}
\multiput(113.67,18.83)(-0.11,-0.20){13}{\line(0,-1){0.20}}
\multiput(112.18,16.24)(-0.12,-0.13){17}{\line(0,-1){0.13}}
\multiput(110.20,14.00)(-0.15,-0.11){16}{\line(-1,0){0.15}}
\multiput(107.82,12.20)(-0.24,-0.12){11}{\line(-1,0){0.24}}
\multiput(105.13,10.90)(-0.41,-0.11){7}{\line(-1,0){0.41}}
\multiput(102.24,10.17)(-1.49,-0.07){2}{\line(-1,0){1.49}}
\multiput(99.25,10.02)(-0.74,0.11){4}{\line(-1,0){0.74}}
\multiput(96.30,10.46)(-0.31,0.11){9}{\line(-1,0){0.31}}
\multiput(93.49,11.49)(-0.20,0.12){13}{\line(-1,0){0.20}}
\multiput(90.94,13.04)(-0.13,0.12){17}{\line(-1,0){0.13}}
\multiput(88.75,15.07)(-0.12,0.16){15}{\line(0,1){0.16}}
\multiput(87.01,17.50)(-0.11,0.25){11}{\line(0,1){0.25}}
\multiput(85.78,20.22)(-0.11,0.49){6}{\line(0,1){0.49}}
\put(85.12,23.13){\line(0,1){2.99}}
\multiput(85.04,26.12)(0.10,0.59){5}{\line(0,1){0.59}}
\multiput(85.56,29.06)(0.11,0.28){10}{\line(0,1){0.28}}
\multiput(86.65,31.84)(0.12,0.18){14}{\line(0,1){0.18}}
\multiput(88.27,34.35)(0.12,0.12){18}{\line(0,1){0.12}}
\multiput(90.36,36.49)(0.16,0.11){15}{\line(1,0){0.16}}
\multiput(92.83,38.17)(0.28,0.12){10}{\line(1,0){0.28}}
\multiput(95.58,39.33)(0.74,0.11){6}{\line(1,0){0.74}}
\put(9.00,36.00){\line(1,-1){22.00}}
\put(31.00,14.00){\line(0,1){0.00}}
\put(31.00,14.00){\line(0,0){0.00}}
\put(31.00,14.00){\line(0,0){0.00}}
\put(31.00,36.00){\line(-1,-1){22.00}}
\put(20.00,40.00){\vector(1,0){1.00}}
\put(35.00,25.00){\vector(0,-1){1.00}}
\put(20.00,10.00){\vector(-1,0){1.00}}
\put(5.00,25.00){\vector(0,1){1.00}}
\put(25.00,20.00){\vector(1,-1){1.00}}
\put(15.00,20.00){\vector(-1,-1){1.00}}
\put(7.00,38.00){\makebox(0,0)[cc]{1}}
\put(33.00,38.00){\makebox(0,0)[cc]{2}}
\put(33.00,12.00){\makebox(0,0)[cc]{3}}
\put(7.00,12.00){\makebox(0,0)[cc]{4}}
\put(49.00,36.00){\line(1,-1){11.00}}
\put(60.00,25.00){\line(1,1){11.00}}
\put(60.00,25.00){\line(0,-1){15.00}}
\put(60.00,40.00){\vector(1,0){1.00}}
\put(45.00,25.00){\vector(0,1){1.00}}
\put(55.00,30.00){\vector(1,-1){1.00}}
\put(65.00,30.00){\vector(-1,-1){1.00}}
\put(47.00,38.00){\makebox(0,0)[cc]{1}}
\put(73.00,38.00){\makebox(0,0)[cc]{2}}
\put(62.00,13.00){\makebox(0,0)[cc]{3}}
\put(58.00,24.00){\makebox(0,0)[cc]{4}}
\put(75.00,25.00){\vector(0,-1){1.00}}
\put(60.00,17.00){\vector(0,1){1.00}}
\put(100.00,25.00){\circle{10.00}}
\put(100.00,30.00){\line(0,1){10.00}}
\put(100.00,20.00){\line(0,-1){10.00}}
\put(100.00,35.00){\vector(0,-1){1.00}}
\put(100.00,15.00){\vector(0,-1){1.00}}
\put(95.00,25.00){\vector(0,-1){1.00}}
\put(105.00,25.00){\vector(0,-1){1.00}}
\put(115.00,25.00){\vector(0,-1){1.00}}
\put(85.00,25.00){\vector(0,1){1.00}}
\put(98.00,37.00){\makebox(0,0)[cc]{1}}
\put(102.00,32.00){\makebox(0,0)[cc]{2}}
\put(98.00,18.00){\makebox(0,0)[cc]{3}}
\put(102.00,13.00){\makebox(0,0)[cc]{4}}
\put(20.00,5.00){\makebox(0,0)[cc]{$\Gamma_1$}}
\put(60.00,5.00){\makebox(0,0)[cc]{$\Gamma_2$}}
\put(100.00,5.00){\makebox(0,0)[cc]{$\Gamma_3$}}
\end{picture}
\caption{}\label{fig-graphsK}
\end{figure}
in fig.\ \ref{fig-dgraphsK} we have instead 
connected decorated graphs for knots of order
2 and degree 1. 
\begin{figure}
\unitlength 1.00mm
\linethickness{0.4pt}
\begin{picture}(100.00,40.00)
\multiput(85.00,40.00)(0.99,-0.10){3}{\line(1,0){0.99}}
\multiput(87.97,39.70)(0.36,-0.11){8}{\line(1,0){0.36}}
\multiput(90.83,38.82)(0.22,-0.12){12}{\line(1,0){0.22}}
\multiput(93.45,37.39)(0.13,-0.11){17}{\line(1,0){0.13}}
\multiput(95.74,35.47)(0.12,-0.15){16}{\line(0,-1){0.15}}
\multiput(97.60,33.14)(0.11,-0.22){12}{\line(0,-1){0.22}}
\multiput(98.96,30.48)(0.12,-0.41){7}{\line(0,-1){0.41}}
\multiput(99.77,27.60)(0.11,-1.49){2}{\line(0,-1){1.49}}
\multiput(100.00,24.63)(-0.09,-0.74){4}{\line(0,-1){0.74}}
\multiput(99.62,21.66)(-0.12,-0.35){8}{\line(0,-1){0.35}}
\multiput(98.67,18.83)(-0.11,-0.20){13}{\line(0,-1){0.20}}
\multiput(97.18,16.24)(-0.12,-0.13){17}{\line(0,-1){0.13}}
\multiput(95.20,14.00)(-0.15,-0.11){16}{\line(-1,0){0.15}}
\multiput(92.82,12.20)(-0.24,-0.12){11}{\line(-1,0){0.24}}
\multiput(90.13,10.90)(-0.41,-0.11){7}{\line(-1,0){0.41}}
\multiput(87.24,10.17)(-1.49,-0.07){2}{\line(-1,0){1.49}}
\multiput(84.25,10.02)(-0.74,0.11){4}{\line(-1,0){0.74}}
\multiput(81.30,10.46)(-0.31,0.11){9}{\line(-1,0){0.31}}
\multiput(78.49,11.49)(-0.20,0.12){13}{\line(-1,0){0.20}}
\multiput(75.94,13.04)(-0.13,0.12){17}{\line(-1,0){0.13}}
\multiput(73.75,15.07)(-0.12,0.16){15}{\line(0,1){0.16}}
\multiput(72.01,17.50)(-0.11,0.25){11}{\line(0,1){0.25}}
\multiput(70.78,20.22)(-0.11,0.49){6}{\line(0,1){0.49}}
\put(70.12,23.13){\line(0,1){2.99}}
\multiput(70.04,26.12)(0.10,0.59){5}{\line(0,1){0.59}}
\multiput(70.56,29.06)(0.11,0.28){10}{\line(0,1){0.28}}
\multiput(71.65,31.84)(0.12,0.18){14}{\line(0,1){0.18}}
\multiput(73.27,34.35)(0.12,0.12){18}{\line(0,1){0.12}}
\multiput(75.36,36.49)(0.16,0.11){15}{\line(1,0){0.16}}
\multiput(77.83,38.17)(0.28,0.12){10}{\line(1,0){0.28}}
\multiput(80.58,39.33)(0.74,0.11){6}{\line(1,0){0.74}}
\put(85.00,34.00){\circle{12.00}}
\put(85.00,28.00){\line(0,-1){18.00}}
\multiput(35.00,40.00)(0.99,-0.10){3}{\line(1,0){0.99}}
\multiput(37.97,39.70)(0.36,-0.11){8}{\line(1,0){0.36}}
\multiput(40.83,38.82)(0.22,-0.12){12}{\line(1,0){0.22}}
\multiput(43.45,37.39)(0.13,-0.11){17}{\line(1,0){0.13}}
\multiput(45.74,35.47)(0.12,-0.15){16}{\line(0,-1){0.15}}
\multiput(47.60,33.14)(0.11,-0.22){12}{\line(0,-1){0.22}}
\multiput(48.96,30.48)(0.12,-0.41){7}{\line(0,-1){0.41}}
\multiput(49.77,27.60)(0.11,-1.49){2}{\line(0,-1){1.49}}
\multiput(50.00,24.63)(-0.09,-0.74){4}{\line(0,-1){0.74}}
\multiput(49.62,21.66)(-0.12,-0.35){8}{\line(0,-1){0.35}}
\multiput(48.67,18.83)(-0.11,-0.20){13}{\line(0,-1){0.20}}
\multiput(47.18,16.24)(-0.12,-0.13){17}{\line(0,-1){0.13}}
\multiput(45.20,14.00)(-0.15,-0.11){16}{\line(-1,0){0.15}}
\multiput(42.82,12.20)(-0.24,-0.12){11}{\line(-1,0){0.24}}
\multiput(40.13,10.90)(-0.41,-0.11){7}{\line(-1,0){0.41}}
\multiput(37.24,10.17)(-1.49,-0.07){2}{\line(-1,0){1.49}}
\multiput(34.25,10.02)(-0.74,0.11){4}{\line(-1,0){0.74}}
\multiput(31.30,10.46)(-0.31,0.11){9}{\line(-1,0){0.31}}
\multiput(28.49,11.49)(-0.20,0.12){13}{\line(-1,0){0.20}}
\multiput(25.94,13.04)(-0.13,0.12){17}{\line(-1,0){0.13}}
\multiput(23.75,15.07)(-0.12,0.16){15}{\line(0,1){0.16}}
\multiput(22.01,17.50)(-0.11,0.25){11}{\line(0,1){0.25}}
\multiput(20.78,20.22)(-0.11,0.49){6}{\line(0,1){0.49}}
\put(20.12,23.13){\line(0,1){2.99}}
\multiput(20.04,26.12)(0.10,0.59){5}{\line(0,1){0.59}}
\multiput(20.56,29.06)(0.11,0.28){10}{\line(0,1){0.28}}
\multiput(21.65,31.84)(0.12,0.18){14}{\line(0,1){0.18}}
\multiput(23.27,34.35)(0.12,0.12){18}{\line(0,1){0.12}}
\multiput(25.36,36.49)(0.16,0.11){15}{\line(1,0){0.16}}
\multiput(27.83,38.17)(0.28,0.12){10}{\line(1,0){0.28}}
\multiput(30.58,39.33)(0.74,0.11){6}{\line(1,0){0.74}}
\put(24.00,14.00){\line(2,5){10.40}}
\put(34.40,40.00){\line(2,-5){10.40}}
\put(50.00,25.00){\vector(0,-1){1.00}}
\put(20.00,25.00){\vector(0,1){1.00}}
\put(35.00,10.00){\vector(-1,0){1.00}}
\put(100.00,25.00){\vector(0,-1){1.00}}
\put(70.00,25.00){\vector(0,1){1.00}}
\put(91.00,34.00){\vector(0,-1){1.00}}
\put(79.00,34.00){\vector(0,-1){1.00}}
\put(85.00,20.00){\vector(0,-1){1.00}}
\put(38.00,38.00){\makebox(0,0)[cc]{1}}
\put(42.00,14.00){\makebox(0,0)[cc]{2}}
\put(28.00,14.00){\makebox(0,0)[cc]{3}}
\put(85.00,38.00){\makebox(0,0)[cc]{1}}
\put(85.00,30.00){\makebox(0,0)[cc]{2}}
\put(87.00,13.00){\makebox(0,0)[cc]{3}}
\put(30.00,30.00){\vector(-1,-2){1.00}}
\put(38.00,30.00){\vector(1,-2){1.00}}
\put(35.00,5.00){\makebox(0,0)[cc]{$\Gamma_1'$}}
\put(85.00,5.00){\makebox(0,0)[cc]{$\Gamma_2'$}}
\end{picture}
\caption{}\label{fig-dgraphsK}
\end{figure}
In all these graphs it is understood
that the outer circle represents the knot.

We will follow the convention of labeling first the external vertices
following the fixed orientation. Only after we have
exhausted the external vertices do we start labeling the internal ones.
In this we have the same freedom as we had before, as well as the
freedom in orienting the internal edges. We divide the algebra
of graphs by the same equivalence relation \eqref{eqrel} we had
before. We keep
calling $\mathcal{D}'$ the quotient.

To a trivalent graph $\Gamma\in\mathcal{D}'$ we can associate the number
\begin{equation}
A_\Gamma(K,M) \doteq \int_{\wcC_{n,t+1}^K} v_0\,
\prod_{(ij)\in\Gamma} \heta_{ij0},
\lab{AKM}
\end{equation}
where $n$ and $t$ are the numbers of external and internal vertices
in $\Gamma$. Notice that by $(ij)$ now we mean only the internal
edges. If we denote by $E$ their number, we have
\begin{equation}
2E=n+3t
\lab{2Ent}
\end{equation}
since the external vertices are univalent as for the internal edges.
This implies that the order of a trivalent graph is 0 and
that the integrand in \eqref{AKM} is actually a
top form on $\wcC_{n,t+1}^K$. Moreover, in this case we have
\begin{equation}
\ord\Gamma = \dfrac{n+t}2.
\lab{ordnt}
\end{equation}

In most cases it is possible to replace $\heta_{ij0}$
by $\heta_{ij}$ in \eqref{AKM}. Actually, we have the following
\begin{Lem}
Unless $(ij)$ 
belongs to an internal loop,
the integral $A_\Gamma(K,M)$ in \eqref{AKM} 
does not change if one replaces
$\heta_{ij0}$ with $\heta_{ij}$.
\lab{lem-loops}
\end{Lem}
\begin{proof}
Since the dimension at the point $0$ is saturated
by the volume form $v_0$, the integration selects
the components of the 2-form $\heta_{ij0}$ that carry
either one form degree at each vertex $i$ and $j$
or two form degrees at one vertex $i$ or $j$ and zero
form degrees on the other. The terms $\heta_{j0}$
and $\heta_{0i}$ contribute to the latter case only.
Therefore, to prove the lemma, it is enough to show
that, unless $(ij)$ belong to an internal loop,
integration selects the component of $\heta_{ij0}$
that carries one form degree on $i$ and one on $j$.

Suppose first that $i$ is an external
vertex. In this case, it is clear that the integral
vanishes by dimensional reasons if we put a zero-
or a 2-form on $i$.

If both $i$ and $j$ are internal we can reason as follows.
Suppose we select the component of $\heta_{ij0}$ that carries
two form degrees on $i$. Call  $j'$ and $j''$ the other
two vertices connected to $i$.
Necessarily, the form on one internal edge, say $\heta_{ij'0}$,
will carry zero form degrees on $i$, while the form on the other
internal edge, say $\heta_{ij''0}$, will carry
one form degree on $i$. Thus, $\heta_{ij'0}$ will carry two form degrees
on $j'$, and so on.
Notice that no vertex can appear twice in this sequence
since otherwise it would carry a 4-form. Moreover, no external
vertex can belong to the sequence since we cannot put a 2-form
on the knot. Thus, this procedure gives a nonvanishing result only
if at some point in this sequence we hit
the vertex $j$; in fact, on $j$ we can put a 2-form
since in $\heta_{ij0}$ we have chosen the component
that carries no form degrees on $j$. 
But this can happen only if $(ij)$ belongs to an
internal loop.
\end{proof}

\begin{Exa}
Referring to the graphs in fig.\ \ref{fig-graphsK}, and taking
into account Lemma \ref{lem-loops}, we have the following integrals:
\[
\begin{array}{ccccc}
A_{\Gamma_1}(K,M) &=& \int_{\wcC_{4,1}^K} v_0\, 
\heta_{130}\,\heta_{240}
&=& \int_{\wcC_{4,0}^K} 
\heta_{13}\,\heta_{24},\\
A_{\Gamma_2}(K,M) &=& \int_{\wcC_{3,2}^K} v_0\, 
\heta_{140}\,\heta_{240}\,\heta_{340}
&=& \int_{\wcC_{3,1}^K} 
\heta_{14}\,\heta_{24}\,\heta_{34},\\
A_{\Gamma_3}(K,M) &=& \int_{\wcC_{2,3}^K} v_0\,
\heta_{120}\,\heta_{230}^2\,\heta_{340}
&=& \int_{\wcC_{2,3}^K} v_0\,
\heta_{12}\,\heta_{230}^2\,\heta_{34}.
\end{array}
\]
\end{Exa}

On $\mathcal{D}'$ we can define a coboundary operator $\delta$
as in the case of decorated graphs for manifolds with the additional
constraint that internal edges connecting external vertices are
not contracted. Thus, $\delta$ contracts external regular edges or
internal regular edges with at least one endpoint internal.
Notice that, if the graph has exactly two external vertices,
there are no external regular edges. We have then an analogue
of Prop.\ \ref{prop-delta}. 

Again we call a cocycle a graph $\Gamma$ killed by $\delta$ and note
that $\delta$ can be restricted to the algebra of connected graph.

\begin{Exa}
The $\Theta$-graph with its outer circle seen as the knot
is a cocycle in the graph cohomology for knots
since it has no regular edges.
\end{Exa}

\begin{Exa}
Referring to the graphs in figs.\ \ref{fig-graphsK} and
\ref{fig-dgraphsK}, we have
\[
\begin{split}
\delta\Gamma_1 &= 4\,\Gamma_1',\\
\delta\Gamma_2 &= 3\,\Gamma_1'+3\,\Gamma_2',\\
\delta\Gamma_3 &= 2\,\Gamma_2'.
\end{split}
\]
Therefore, the combination
\[
\Gamma=
\frac14\,\Gamma_1 - \frac13\,\Gamma_2+\frac12\,\Gamma_3
\]
is a cocycle in the graph cohomology for knots.

Notice that $\delta\Gamma_1$ differs here from
what we obtained in Example \ref{coGamma} since there we had
to contract also the edges $(13)$ and $(24)$. \lab{exa-coGammak}
\end{Exa}

Finally we can state the following
\begin{Thm}
If $K$ is a knot in the rational homology 3-sphere $M$ and\/
$\Gamma$ a connected, trivalent cocycle 
in the graph cohomology for
knots, then there exists a constant $\mu(\Gamma)$ such that the
quantity
\[
I_\Gamma(K,M) = A_\Gamma(K,M) + \mu(\Gamma)\, \sln(K,M)
\]
is a knot invariant. Moreover, $\mu(\Gamma)=0$ if\/ $\ord\Gamma$ is
even.
\lab{thm-GammaK}
\end{Thm}
For $M=\bbR^3$, the analogous theorem was proved in \cite{BT};
the vanishing of $\mu(\Gamma)$ for $\ord\Gamma$ even in this case
was proved in \cite{AF}. The simplest knot invariant---which
corresponds to the cocycle described in Example 
\ref{exa-coGammak}---had previously been described in
\cite{GMM} and \cite{BN}.

\begin{Rem}
In \cite{BN} a product structure is introduced on the algebra
of graphs. With this product,
\[
A_{\Gamma_1\cdot\Gamma_2}(K,M) = A_{\Gamma_1}(K,M)\,
A_{\Gamma_2}(K,M).
\]
This explains why it is enough to restrict our attention
to connected graphs only. (Actually, in \cite{AF} it
is shown that as for the computation of $\mu(\Gamma)$
it is enough to consider ``primitive" graphs, namely,
decorated graphs for knots such that their underlying
graphs are connected after removing all external edges.)
\end{Rem}

\subsection{Proof of Thm.\ \ref{thm-GammaK}}
Since the integrand form is closed by construction, then
$dI_\Gamma$---as a 1-form on the unit interval $I$---will
be given just in terms of boundary integrals.

We first consider the faces corresponding to the collapse of vertices
(necessarily internal or carrying the volume form)
at an internal vertex. In this case we can use the same arguments
we used in subsection \ref{ssec-ptGamma} to prove Thm.\ \ref{thm-Gamma}.
Essentially we come to the same conclusions of Lemma \ref{lem-Gamma},
with one important difference: since there are always at least two
external vertices and the diagram is connected, it can never happen that
all points come together at an internal vertex. Therefore, we are left
only with case 1 of Lemma \ref{lem-Gamma}, and this is taken care of
by the action of the coboundary operator $\delta$.

Now consider a collapse at an external vertex. Here 
both internal and external vertices can come together. 
Notice that, if the point
0 is involved, the form vanishes since $v$ is a 3-form and
$S^1\times I$ is 2-dimensional; so we must only consider
the case when $q$ external and $s$ internal vertices (with $q\ge1$,
$s\ge0$ and $q+s\ge2$) come together
at a point $x$ on the knot. Let us denote by 
$\mathcal{S}$ the corresponding face. This is a bundle
over $I$ whose fiber is a component of $\de \wcC_{n,t+1}^{K_\tau}$. 
However,
we can also think of $\mathcal{S}$ as
a bundle over $\hwcC_{n-q+1,t-s+1}^K$
whose projection we denote by $\pi^\de$. Its $(q+3s-2)$-dimensional
fiber can be described as given by $q$ copies of $T_xS^1$ 
and $s$ copies of $T_{K_\tau(x)}M$ up to global translations along $T_xS^1$
and scalings, and with the diagonals blown up. Moreover, since we
are considering only one connected component of the configuration
spaces, we must fix an ordering of the $q$ points on $T_xS^1$.

To give an explicit description of this fiber, we consider
the following commutative diagram:
\[
\begin{CD}
\mathcal{S} @>{\hat f}>> \pi^{-1}P\times_{SO(3)}(S^2\ltimes F)\\
@V{\pi^\de}VV   @VV{\hat\pi^\de}V\\
\hwcC_{n-q+1,t-s+1}^K @>>{f}> \pi^{-1}P\times_{SO(3)}S^2
\end{CD}
\]
Here $\pi$ is the projection $M\times I\to M$, and $f$ is the 
composition of the projection 
\[
\pi_x:\hwcC_{n-q+1,t-s+1}^K\to S^1\times I
\]
to the point
where the collapse has happened with the map 
\[
\hat\psi_K : \begin{array}[t]{ccc}
S^1\times I &\to& \pi^{-1}P\times_{SO(3)} S^2,\\
(\alpha,\tau) &\mapsto& (\psi_{K_\tau},\tau)
\end{array}
\]
with $\psi_K$ defined in \eqref{defpsi}. (Remember that 
$\de C_2(M)=P\times_{SO(3)}S^2$.)

The space $S^2\ltimes F$ is defined as follows.
To each $\ba\in S^2$ we associate the imbedding
\[
\ba : \begin{array}[t]{ccc}
\bbR &\to&\bbR^3,\\
x &\mapsto &\ba\,x.
\end{array}
\]
Then we call $F^\ba$ the configuration space $\wcC_{q,s}^\ba$
modulo translations and scalings; that is,
if we denote by $x_1,\dots, x_q$
the $q$ coordinates in $\bbR$ (with $x_1<x_2<\dots<x_q$) and
by $\bx_{q+1},\dots,\bx_{q+s}$  the $s$ coordinates in $\bbR^3$,
we divide $\wcC_{q,s}^\ba$ by the translations
\begin{align*}
x_i &\to x_i+\xi,\\
\bx_j &\to \bx_j + \xi\ba,\\
\intertext{with $\xi\in\bbR$, and by the scalings}
x_i &\to \lambda x_i,\\
\bx_j &\to \lambda \bx_j,
\end{align*}
with $\lambda\in\bbR^*$. By $S^2\ltimes F$ we then mean
the pairs $(\ba,F^\ba)$ with $\ba\in S^2$.

The action of $SO(3)$ on $S^2\ltimes F$
is just the defining action on the copies of $\bbR^3$ and on
$S^2$, and the trivial action on the copies of $\bbR$.

Next consider the form $\lambda\in\Omega^*(\mathcal{S})$ given
by the restriction to this face of one summand of
the integrand form in $I_\Gamma$.
We can split $\lambda$ as
\[
\lambda = \lambda_1\,\pi^{\de*}\lambda_2.
\]
Then we have the following
\begin{Lem}
If $\lambda_1$ is a vertical form on $\mathcal{S}$
given by the restriction to this face of a sum of products
of forms $\heta_{ij}$, then $\pi^\de_*\lambda_1$ vanishes
unless it is a zero-form, in which case it is a constant, or
is a 2-form, in which case it is a multiple of $f^*\eta$.
\lab{lem-f*eta}
\end{Lem}
\begin{proof}
We can write $\lambda_1=\hat f^*\hat\lambda_1$, where $\hat\lambda_1$
is a sum of products of pullbacks of the form 
$\eta\in\pi^{-1}P\times_{SO(3)}S^2$
through the projections $\pi_{ij}$ we are going to describe. First
define
\[
\pi_{ij}:
S^2\ltimes F \to S^2,
\]
by
\[
\pi_{ij}(\ba,x_1,\dots,x_q,\bx_{q+1},\dots,\bx_{q+s}) =
\begin{cases}
\ba\sgn(j-i) & \text{if $i,j\le q$},\\
\rho(\bx_j-\ba x_i) & \text{if $i\le q, j>q$},\\
\rho(\ba x_j-\bx_i) & \text{if $i > q, j\le q$},\\
\rho(\bx_j-\bx_i) & \text{if $i,j>q$},
\end{cases}
\]
with
\[
\rho(\bx) = \dfrac\bx{|\bx|}.
\]
Since the $\pi_{ij}$'s are equivariant they descend to
\[
\pi_{ij} : \pi^{-1}P\times_{SO(3)}(S^2\ltimes F)
\to \pi^{-1}P\times_{SO(3)} S^2.
\]
Now consider the commutative diagram:
\[
\begin{CD}
\pi^{-1}P\times_{SO(3)}(S^2\ltimes F) @<{\Bar{\Hat p}}<<
\pi^{-1}P\times (S^2\ltimes F) @>{\pi_{ij}}>> \pi^{-1}P\times S^2\\
@V{\hat\pi^\de}VV     @VV{\Bar{\Hat \pi}^\de}V \\
\pi^{-1}P\times_{SO(3)} S^2 @<<{\hat p}< \pi^{-1}P\times S^2
\end{CD}
\]
The form $\Bar{\Hat p}^*\hat\lambda_1$ is given
in terms of pullbacks  of the form 
$\bar\eta=p^*\eta$, defined in \eqref{bpeta}, 
through the projections $\pi_{ij}$. Since the
maps $\pi_{ij}$'s and $\Bar{\Hat \pi}^\de$ act as the identity
on $\pi^{-1}P$, we conclude that 
$\Bar{\Hat \pi}^\de_*\Bar{\Hat p}^*\hat\lambda_1$ is a
polynomial in $\theta$ and $d\theta$ with coefficients in
$\Omega^*(S^2)$. Moreover, we know that it is basic. Therefore,
it must be either a constant or a multiple of the form $\bar\eta$.
As a consequence, $\hat\pi^\de_*\hat\lambda_1$ is either a constant
or a multiple of the form $\eta$.
\end{proof}

Now we compute the degree of $\pi^\de_*\lambda_1$.
If $i$ and $j$ are in $\mathcal{S}$ then $\heta_{ij0}$
reduces to $\heta_{ij}$. If at least one vertex $i$ or $j$
is not in $\mathcal{S}$, then $\heta_{ij0}$ is basic in $\mathcal{S}$
(that is, it contributes to $\lambda_2$). Thus, the degree of $\lambda_1$
is equal to the number of internal edges connecting vertices inside $\mathcal{S}$.
Let us denote this number by $e$, and let $e_0$ be the number of
internal edges connecting a point in $\mathcal{S}$ with a point outside.
Then we have
\[
2e+e_0 = q+3s,
\]
since in $\mathcal{S}$ we have $q$ univalent and $s$ trivalent
vertices (as for the internal edges).
Since the fiber dimension is $q+3s-2$ we conclude that
\[
\deg\pi^\de_*\lambda_1=2-e_0.
\]

By Lemma \ref{lem-f*eta}, we then see that $\pi^\de_*\lambda_1$
vanishes unless $e_0=0$ or $e_0=2$. In the former instance
all points but the point 0 collapse on the knot. The point 0
is now completely disconnected from the point on the knot
where all points have collapsed; therefore,
we can integrate the volume form on 0. We are then left with
a multiple of $f^*\eta=\psi_{K_\tau}^*\eta$ to be integrated over 
$\wcC_{1,0}^K=S^1$. But this can also be written as a multiple
of $d\sln(K_\tau,M)$ by \eqref{dsln}.

To prove that we do not have this contribution if $\ord\Gamma$ is
even, we consider the involution $\phi$ that acts: as
the antipodal map on $S^2$, trivially on the $q=n$ copies of
$\bbR$, and as the reflection with respect to the origin on
the $s=t$ copies of $\bbR^3$. Since the maps $\pi_{ij}$ are equivariant
and $\phi^*\eta=-\eta$ by Prop.\ \ref{thm-eta}, we have
\[
\phi^*\hat\lambda_1 = (-1)^E\,\hat\lambda_1,
\]
where $E$ is the number of internal edges. The map $\hat\pi^\de$
is also equivariant; however, the orientation of the fiber is reversed
if $t$ is odd, so
\[
\phi^*\hat\pi^\de_*\hat\lambda_1 = (-1)^{E+t}\,\hat\pi^\de_*
\hat\lambda_1.
\]
On the other hand
we know that $\hat\pi^\de_*\hat\lambda_1$ is proportional to
$\eta$, so
\[
\phi^*\hat\pi^\de_*\hat\lambda_1 = -\hat\pi^\de_*
\hat\lambda_1.
\]
This means that $\hat\pi^\de_*\hat\lambda_1$ vanishes if
$E+t$ is even. By \eqref{2Ent} and \eqref{ordnt}, this is
equivalent to $\ord\Gamma$ even.

To complete our proof, we must finally consider the case $e_0=2$.
In this case the fiber dimension is equal to the degree of $\lambda_1$.
As in the proof of Lemma \ref{lem-f*eta}, we write
$\lambda_1=f^*\hat\lambda_1$ and see that $\hat\pi^\de_*$
now selects the part of degree 0 in $\theta$. Therefore, we can
use the same arguments used in \cite{BT} to prove that
$\pi^\de_*\lambda_1$ vanishes unless $q+s=2$. This
case is taken care of by the coboundary operator $\delta$. (Notice that
if $\Gamma$ has exactly two external vertices, say 1 and 2, 
and we are considering their collapse, then we get two opposite
contributions as 1 approaches 2 from one or the other side, provided
that 1 and 2 are not connected by one internal edge.)
 
This concludes the proof of Thm.\ \ref{thm-GammaK}.

\end{document}